\documentclass[12pt]{amsart}

\usepackage{graphicx,verbatim
 }
\usepackage{amsmath} \usepackage{amsfonts}

\makeindex

\newcommand{\ECH}{\mathrm{\mathop{conv_{\mathcal{E}}}}}
\newcommand{\HCH}{\mathrm{\mathop{conv_{\mathcal H}}}}

\newcommand{\mc}{\mathcal}
\newcommand{\mrm}{\mathrm}

\newcommand{\ol}{\overline}
\newcommand{\cl}{\overline}

\newcommand{\KP}{\mathcal{KP}}
\newcommand{\kp}{\mathcal{KP}}
\newcommand{\KPP}{\mathcal{KPP}}

\newcommand{\complex}{\mathbb{C}}
\newcommand{\Complex}{\mathbb{C}}
\newcommand{\rsphere}{\mathbb{C}^\infty}
\newcommand{\sphere}{\Complex^\infty}
\newcommand{\ucirc}{\mathbb{S}^1}

\newcommand{\B}{\mathfrak{B}}

\newcommand{\real}{\mathbb{R}}

\newcommand{\zed}{\mathbb{Z}}
\newcommand{\Disk}{\mathbb{D}}
\newcommand{\disk}{\Delta}
\newcommand{\wh}{\widehat}
\newcommand{\e}{\varepsilon}
\newcommand{\0}{\emptyset}

\newcommand{\var}{\mathrm{var}}
\newcommand{\ind}{\mathrm{ind}}
\newcommand{\dg}{\mathrm{degree}}
\newcommand{\win}{\mathrm{win}}
\newcommand{\bd}{\partial}
\newcommand{\im}{\mathrm{Im}}
\newcommand{\re}{\mathrm{Re}}
\newcommand{\pr}{\mathrm{Pr}}
\newcommand{\dm}{\mathrm{diam}}

\newcommand{\sm}{\setminus}
\newcommand{\degree}{\text{degree}}
\newcommand{\set}[1]{\left\{#1\right\}}
\newcommand{\hX}{\hat{X}}
\newcommand{\Int}{\mathrm{Int}}
\newcommand{\Sh}{\mathrm{Sh}}
\newcommand{\al}{\alpha}
\newcommand{\tf}{\tilde{f}}
\newcommand{\tU}{U^\infty}
\newcommand{\rg}{\mathrm{\mathbf{g}}}
\newcommand{\rh}{\mathrm{\mathbf{h}}}
\newcommand{\fg}{\mathfrak{g}}

\newcommand{\rG}{\mathrm{G}}
\newcommand{\rH}{\mathrm{H}}
\newcommand{\hx}{\hat{x}}
\newcommand{\vp}{\varphi}
\newcommand{\fS}{\mathfrak{S}}
\newcommand{\ta}{\tilde{a}}
\newcommand{\tb}{\tilde{b}}
\newcommand{\tB}{\tilde{B}}

\newtheorem{thm}{Theorem}[section]
\newtheorem{lem}[thm]{Lemma}
\newtheorem{cor}[thm]{Corollary}
\newtheorem{prop}[thm]{Proposition}
\newtheorem{defn}[thm]{Definition}
\newtheorem{rem}[thm]{Remark}

\begin{document}
\title[Fixed points in planar continua]{The plane fixed point problem}
\dedicatory{Dedicated to Harold Bell}
\author[Fokkink]{Robbert~Fokkink }
\address{Delft University,
Faculty of Information and Systems, P.O. Box 5031, 2600 GA Delft,
Netherlands } \email[Robbert Fokkink]{R.J.Fokkink@its.tudelft.nl}
\author[Mayer]{John C.~Mayer}
\address[John C.~Mayer and Lex G.~Oversteegen]
{Department of Mathematics\\ University of Alabama at Birmingham\\
Birmingham, AL 35294-1170} \email[John
C.~Mayer]{mayer@math.uab.edu}
\author[Oversteegen]{Lex G.~Oversteegen}
\email[Lex G.~Oversteegen]{overstee@math.uab.edu}
\author[Tymchatyn]{E.~D.~Tymchatyn}

\address[E.~D.~Tymchatyn]{Department of Mathematics and Statistics\\
University of Saskatchewan\\
Saskatoon, Saskatchewan, Canada S7N 0W0}
\email[E.~D.~Tymchatyn]{tymchat@math.usask.ca}
\thanks{The third named author was supported in part by grant
NSF-DMS-0405774 and the last named author by NSERC 0GP0005616.}
\date{\today}

\begin{abstract}
In this paper we  present  proofs of basic results,
including those developed
so far  by H.~Bell, for the plane fixed point problem. Some of these results had been
announced much earlier by Bell but without accessible  proofs. We define the concept of the
variation of a map on a simple closed curve and relate it to the index of the map on
that curve: Index = Variation + 1. We develop a prime end theory through hyperbolic
chords in maximal round  balls contained in the complement of a non-separating plane continuum $X$.
We define the concept of an {\em outchannel} for a fixed point free map which carries
the boundary of $X$ minimally  into itself and  prove  that such a map has a \emph{unique}
outchannel, and that outchannel must have variation $=-1$. We also extend Bell's linchpin theorem
for a foliation of a simply connected domain, by closed convex subsets, to arbitrary
domains in the sphere.

We introduce the notion of an oriented map of the plane. We show that the perfect oriented maps
of the plane coincide with confluent (that is composition of monotone and open) perfect maps
of the plane. We obtain a fixed point theorem for  positively oriented, perfect maps
of the plane. This generalizes results announced by Bell in 1982 (see also \cite{akis99}).
It follows that if $X$ is invariant  under an oriented  map $f$, then $f$ has a point of period
at most two in $X$.

\end{abstract}

\maketitle
\tableofcontents

\section{Introduction}

We denote the plane  by $\complex$, \index{complex@$\complex$} the Riemann
sphere by  $\rsphere=\complex\cup\{\infty\}$, \index{complex@$\rsphere$} the real line by
$\real$ \index{real@$\real$} and the unit circle by $S^1=\real/\zed$. Let $X$ be a plane continuum.
Since $\complex$ is locally connected and $X$ is closed, complementary domains of $X$ are open.
By $T(X)$ \index{TX@$T(X)$} we denote the {\em
topological hull} \index{topological hull} of $X$ consisting of $X$ union all of its
bounded complementary domains.  Thus, $U^\infty=\rsphere\sm T(X)$
\index{U@$U^\infty$} is a
simply-connected open domain containing $\infty$.  The following is a
long-standing question in topology.

{\bf Fixed Point Question}: {\em ``Does a continuous function
taking a non-separating plane continuum into itself always have a
fixed point?"}

It is easy to see that a map of a plane continuum to itself can be extended to a perfect map of the plane.
We study  the slightly more general question, ``Is there a
plane continuum $Z$ and a  perfect continuous function
$f:\complex\to\complex$ taking $Z$ into $T(Z)$ with no fixed
points in $T(Z)$?" A Zorn's Lemma argument shows that if one
assumes the answer is ``yes," then there is a subcontinuum
$X\subset Z$ minimal with respect to these properties.
It will follow from Theorem~\ref{densechannel} that for such a
minimal continuum, $f(X)=X=\bd T(X)$ (though it may not be the case
that $f(T(X))\subset T(X)$).  Here $\partial T(X)$ denotes the boundary of $T(X)$\index{del@$\partial$ boundary operator}.
 We recover Bell's result \cite{bell67} (see also Sieklucki
\cite{siek68}, and Iliadis \cite{ilia70}) that  the
boundary of $X$ is  indecomposable (with a dense channel,
explained later).

In this paper we use tools first developed by Bell to elucidate the action of a fixed
point free map (should one exist). We are indebted to Bell for sharing his insights  with us.
Many of the results of this paper were first obtained by him. Unfortunately, many of the proofs were not accessible.
We believe that they deserve to be developed in order to be
useful to the mathematical community. The results of this paper are also crucial to several recent results
regarding the extension of isotopies of plane continua \cite{overtymc07}, the existence of fixed points
for branched covering maps of the plane \cite{blokoverto1}, fixed points in non-invariant plane continua
\cite{blokover08a}, the existence of locally connected models
for all connected Julia sets of complex polynomials \cite{blokover08b} and an estimate on the number of
attracting and neutral periodic orbits of complex polynomials \cite{blokoverto2}.

 We have stated many of
these results using existing notions such as prime ends.
We introduce Bell's notion of variation and prove his theorem that index equals variation +1;
Theorem~\ref{I=V+1}. We also extended Bell's linchpin
Theorem~\ref{Hypmain} for simply connected domains  to arbitrary domains in the sphere and given a proof using
 an elegant argument due to Kulkarni and Pinkall \cite{kulkpink94}.  Our version of this  theorem
(Theorem~\ref{KPthrm}) is
essential for the results later in the paper.     Theorem~\ref{outchannel} (Unique Outchannel)
is a new result due to Bell.  Complete proofs of Theorems~\ref{I=V+1}, \ref{lollipop},
 \ref{Hypmain}
and \ref{outchannel}
appear in print for the first time.

The classical fixed point question asks whether each map of a
non-separating plane continuum into itself must have a fixed
point. Cartwright and Littlewood \cite{cartlitt51} showed that the
answer is yes if the map can be extended to an
orientation-preserving homeomorphism of the plane. It was 25
years before Bell \cite{bell78} extended this to the class of all
homeomorphisms of the plane.  Bell
announced in 1984 (see also Akis \cite{akis99}) that the
Cartwright-Littlewood Theorem can be extended to the class of all
holomorphic maps of the plane. These maps behave like
orientation-preserving homeomorphisms in the sense that they
preserve local orientation.  Compositions of
open, perfect  and of monotone, perfect surjections of the plane  are confluent and naturally
decompose into two classes, one of which preserves and the other
of which reverses local orientation. We show that any confluent map of the plane
 is itself a composition of a monotone and a light-open map of the plane.
 We  also show that  an oriented map of the plane induces a map to the
circle of prime ends of an acyclic continuum   from the
circle of prime ends of a component of its pre-image. Finally we will show that each
invariant non-separating plane continuum, under a
positively-oriented map of the plane, must contain a fixed point. It follows that
any confluent map of the plane has a point of period at most two in any non-separating
invariant sub-continuum.

For the convenience of the reader we have included an index at the end of the paper.

\section{Tools}
Let
$p:\real\to S^1$ denote the covering map $p(x)=e^{2\pi ix}$. Let
$g:S^1\to S^1$ be a map. By the \emph{degree} \index{degree} of the map $g$,
\index{degree@$\dg(g)$}
denoted by $\dg(g)$, we mean the number $\hat{g}(1)-\hat{g}(0)$,
where $\hat{g}:\real\to\real$ is a lift of the map $g$ to the
universal covering space $\real$ of $S^1$ (i.e.,
$p\circ\hat{g}=g\circ p$).
 It is well-known that $\dg(g)$ is
independent of the choice of the lift.

\subsection{Index}\label{defindex}

Let $g:S^1\to\complex$ be a map and
$f:g(S^1)\to\complex$ a fixed point free map.   Define
the map $v:S^1\to S^1$ by
$$v(t)=\frac{f(g(t))-g(t)}{|f(g(t))-g(t)|}.$$

 Then the map $v:S^1\to S^1$  lifts to a map $\wh
v:\real\to\real$. Define the
{\em index of $f$ \index{index} with respect to $g$}, denoted
\index{index1@$\ind(f,g)$} $\ind(f,g)$ by
$$\ind(f,g)=\wh v(1)-\wh v(0)=\dg(v).$$

Note that $\ind(f,g)$ measures the net number of revolutions of the vector\linebreak
$f(g(t))-g(t)$ as $t$ travels through the unit circle one revolution in the positive direction.

\begin{rem}\label{remark} (a) If $g:S^1 \to\complex$
is a constant map with $g(S^1)=c$ and $f(c)\ne c$, then $\ind(f,g)=0$. \\
(b)If $f$ is a constant map and $f(\complex)=w$ with $w\not\in g(S^1)$, then $\ind(f,g)=\win(g, S^1,w)$,
the
\emph{ winding number of $g$ about $w$}\index{win@$\win(g,S^1,w)$}.
In particular, if $f:S^1\to T(S^1)\setminus S^1$ is a constant map, then $\ind(f,id|_{S^1})= 1$,
where $id|_{S^1}$ is the identity map on $S^1$\index{id@$id$ identity map}.
\end{rem}

Suppose $S\subset\Complex$ is a simple closed curve and $A\subset S$ is a subarc
\index{order!on subarc of simple closed curve}
of $S$ with endpoints $a$ and $b$.
Then we write $A=[a,b]$ if $A$ is the arc obtained by traveling in the counter-clockwise direction
from the point $a$ to the point $b$ along $S$. In this case we denote by $<$ the linear order
on the arc $A$ such that $a<b$.  \index{counterclockwise order!on an arc in a simple closed curve}
We  will call the order $<$ the \emph{counterclockwise
order on}   $A$. Note that $[a,b]\ne[b,a]$.

More generally, for any arc $A=[a,b]\subset S^1$, with $a<b$ in the counterclockwise order,
define the {\em fractional index} \index{index!fractional} \cite{brow90}
\index{index2@$\ind(f,g|_{[a,b]})$} of $f$ on the
sub-path $g|_{[a,b]}$  by $$\ind(f,g|_{[a,b]})=\wh v(b)-\wh
v(a).$$ While, necessarily, the index of $f$ with respect to $g$ is
an integer, the fractional index of $f$ on $g|_{[a,b]}$ need not
be. We shall have occasion to use fractional index in the proof of
Theorem~\ref{I=V+1}.

\begin{prop} \label{fracindex} Let $g:S^1\to\complex$ be a map with $g(S^1)=S$, and
suppose $f:S\to\complex$ has no fixed points on $S$. Let
$a\not=b\in S^1$ with $[a,b]$ denoting the counterclockwise subarc
on $S^1$ from $a$ to $b$  (so $S^1=[a,b]\cup[b,a]$). Then
$\ind(f,g)=\ind(f,g|_{[a,b]})+\ind(f,g|_{[b,a]})$.
\end{prop}

\subsection{Stability of Index} \label{compofind}
The following standard theorems and observations about the
stability of index under a fixed point free homotopy are
consequences of the fact that index is continuous and
integer-valued.

\begin{thm}\label{fpfhomotopy}  Let $h_t:S^1\to\complex$ be a homotopy.
If $f:\cup_{t\in[0,1]} h_t(S^1)\to \complex$ is fixed point free, then
 $\ind(f,h_0)=\ind(f,h_1)$.
\end{thm}

An embedding $g:S^1\to S\subset\complex$ is \emph{orientation
preserving} \index{embedding!orientation preserving} \index{orientation preserving!embedding}
if $g$ is isotopic to the identity map $id|_{S^{1}}$.
It follows from Theorem~\ref{fpfhomotopy} that if $g_1,g_2:S^1\to S$
are orientation preserving homeomorphisms and $f:S\to\complex$ is a
fixed point free map, then $\ind(f,g_1)=\ind(f,g_2)$. Hence we can denote
$\ind(f,g_1)$ by $\ind(f,S)$ \index{index3@$\ind(f,S)$} and if $[a,b]$ is a positively oriented subarc
of $S^1$ we denote $\ind(f,g_1|_{[a,b]})$ by $\ind(f,g_1([a,b]))$, \index{index0@$\ind(f,A)$}
by some abuse of notation when the extension of $g_1$  over $S^1$ is
understood.

\begin{thm} \label{fpfhomotopy2} Suppose $g:S^1\to\complex$ is a map
with $g(S^1)=S$, and $f_1,f_2:S\to\complex$ are homotopic maps
such that each level of the homotopy is fixed point free on $S$.
Then $\ind(f_1,g)=\ind(f_2,g)$.
\end{thm}

In particular, if $S$ is a simple closed curve and
$f_1,f_2:S\to\complex$ are  maps such that there is a homotopy
$h_t:S\to\complex$ from $f_1$ to $f_2$ with $h_t$ fixed point free
on $S$ for each $t\in[0,1]$, then $\ind(f_1,S)=\ind(f_2,S)$.

\begin{cor} \label{mapinhull} Suppose $g:S^1\to\complex$ is
an orientation preserving embedding with $g(S^1)=S$, and
$f:S\to T(S)$ is a fixed point free map.  Then $\ind(f,g)=\ind(f,S)=1$.
\end{cor}

\begin{proof}  Since $f(S)\subset T(S)$ which is a disk with boundary $S$ and
$f$ has no fixed point on $S$, there is a  fixed point free
homotopy of $f|_S$ to a constant map $c:S\to \complex$ taking $S$
to a point in $T(S)\setminus S$. By Theorem~\ref{fpfhomotopy2},
$\ind(f,g)=\ind(c,g)$.  Since $g$ is orientation preserving it follows from
Remark~\ref{remark} (b) that $\ind(c,g)= 1$.
\end{proof}

\begin{thm}  \label{fpthm} Suppose $g:S^1\to\complex$ is a map
with $g(S^1)=S$, and $f:T(S)\to\complex$ is a map such that
$\ind(f,g)\not= 0$, then $f$ has a fixed point in $T(S)$.
\end{thm}

\begin{proof}  Notice that $T(S)$ is a locally connected, non-separating, plane continuum
 and, hence, contractible.  Suppose $f$ has no fixed point in $T(S)$.  Choose
point $q\in T(S)$.  Let $c:S^1\to \complex$ be the constant map
$c(S^1)=\{q\}$.  Let $H$ be a homotopy from $g$ to $c$ with image in $T(S)$.
Since $H$ misses the fixed point set of $f$,
Theorem~\ref{fpfhomotopy} and Remark~\ref{remark} (a) imply  $\ind(f,g)=\ind(f,c)=0$.
\end{proof}

\subsection{Variation} \label{compofvar}In this subsection we introduce the notion
of variation of a map on an arc and relate it to winding number.

\begin{defn} [Junctions]\label{junction}\index{junction}
The {\em standard junction} $J_O$ is the union of the three rays
$J^i_O=\{z\in\complex\mid z=re^{i\pi/2},\ r\in[0,\infty)\}$,
$J^+_O=\{z\in\complex\mid z=r,\ r\in[0,\infty)\}$,
$J^-_O=\{z\in\complex\mid z=re^{i\pi},\ r\in[0,\infty)\}$, having
the origin $O$ in common.   A {\em junction}
$J_v$ is the image of $J_O$ under any orientation-preserving
homeomorphism $h:\complex\to\complex$ where $v=h(O)$.
We will often suppress $h$ and refer to $h(J^i_O)$ as $J^i_v$, and
similarly for the remaining rays in $J_v$.  Moreover, we require that for each neighborhood
$W$ of $v$, $d(J^+_v\sm W, J^i_v\sm W)>0$.
 \end{defn}

\begin{defn}[Variation on an arc] \label{vararc} Let $S\subset\Complex$ be a simple closed curve,
$f:S\to\complex$ a map
and $A=[a,b]$ a subarc of $S$ such that $f(a),f(b)\in T(S)$ and
$f(A)\cap A=\0$.  We define the {\em variation of $f$ on $A$ with
respect to $S$}\index{variation!on an arc}, denoted $\var(f,A,S)$, by the following
algorithm:
\begin{enumerate}
\item Let $v\in A$ and let $J_v$ be a junction with $J_v\cap S=\{v\}$.
\item\label{crossings} {\em Counting crossings:} Consider the set
$M=f^{-1}(J_v)\cap [a,b]$. Each time a point of
$f^{-1}(J^+_v)\cap [a,b]$ is immediately followed in $M$, in the
counterclockwise order $<$ on $[a,b]\subset S$, by a point of $f^{-1}(J^i_v)$ count
$+1$ and each time a point of $f^{-1}(J^i_v)\cap [a,b]$ is
immediately followed in $M$  by a
point of $f^{-1}(J^+_v)$ count $-1$. Count no other crossings.
\item The sum of the crossings found above is the variation    $\var(f,A,S)$.  \index{variation@$\var(f,A,S)$}
\end{enumerate}
\end{defn}

  Note that $f^{-1}(J^+_v)\cap [a,b]$ and $f^{-1}(J^i_v)\cap [a,b]$
  are disjoint closed sets in $[a,b]$. Hence, in (\ref{crossings}) in
  the above  definition, we count only a finite number of crossings
  and $\var(f,A,S)$ is an integer. Of course, if $f(A)$ does not meet both
  $J^+_v$ and $J^i_v$, then $\var(f,A,S)=0$.

 If $\al:S\to\complex$ is any map
such that $\al|_A=f|_A$ and $\al(S\setminus (a,b))\cap J_v=\0$,
then $\var(f,A,S)=\win(\al,S,v)$. In particular, this condition is satisfied if
$\al(S\sm (a,b))\subset T(S)\setminus\{v\}$.
The
invariance of winding number under suitable homotopies implies
that the variation $\var(f,A,S)$ also remains invariant under such
homotopies. That is, even though the specific crossings in
(\ref{crossings}) in the algorithm may change, the sum remains
invariant. We will state the  required results about variation
below without proof. Proofs can also be obtained directly by using
the fact that $\var(f,A,S)$ is integer-valued and continuous under
suitable homotopies.

\begin{prop} [Junction Straightening] \label{straightjunction} \label{int}
Let $S\subset\Complex$ be a simple closed curve,
$f:S\to\complex$ a map
and $A=[a,b]$ a subarc of $S$ such that $f(a),f(b)\in T(S)$ and
$f(A)\cap A=\0$.  Any
two junctions $J_v$ and $J_u$  with $u,v\in A$ and $J_w\cap S=\{w\}$
 for $w\in\{u,v\}$
 give the same value for  $\var(f,A,S)$.  Hence $\var(f,A,S)$ is independent
 of the particular junction used  in Definition~\ref{vararc}.
 \end{prop}

The
computation of $\var(f,A,S)$ depends only upon the crossings of
the junction $J_v$  coming from a proper compact subarc of the open arc
$(a,b)$. Consequently, $\var(f,A,S)$ remains invariant under
homotopies $h_t$ of $f|_{[a,b]}$ in the complement of $\{v\}$ such that  $h_t(a),h_t(b) \not\in J_v$
 for all $t$. Moreover, the
computation is stable under an isotopy $h_t$ of the plane that moves
the entire junction $J_v$ (even off $A$), provided in the
isotopy $h_t(v)\not\in f(A)$ and $f(a),f(b)\not\in h_t(J_v)$ for all $t$.

In case $A$ is an open arc $(a,b)\subset S$ such that
$\var(f,\cl{A},S)$ is defined, it will be convenient to denote
$\var(f,\cl{A},S)$ by $\var(f,A,S)$

The following Lemma follows immediately from the definition.
\begin{lem} \label{summ} Let $S\subset\Complex$ be a simple closed curve. Suppose that $a<c<b$ are three points in $S$ such that $\{f(a),f(b),f(c)\}\subset T(S)$ and $f([a,b])\cap [a,b]=\0$.
Then $\var(f,[a,b],S)=\var(f,[a,c],S)+\var(f,[c,b],S)$.
\end{lem}

\begin{defn}[Variation on a finite union of arcs] \label{partition}
\index{variation!on finite union of arcs}
Let $S\subset\Complex$ be a simple closed curve and $A=[a,b]$ a subcontinuum of
$S$ with partition a finite set $F=\{a=a_0<a_1<\dots<a_n=b\}$.
For each $i$ let $A_i=[a_i,a_{i+1}]$.
 Suppose that $f$
satisfies $f(a_i)\in T(S)$ and $f(A_i)\cap A_i=\0$ for each $i$.
We define the \emph{variation of $f$ on $A$ with respect to  $S$},
denoted $\var(f,A,S)$, by
$$\var(f,A,S)=\sum_{i=0}^{n-1} \var(f,[a_i,a_{i+1}],S).$$  In
particular, we include the possibility that  $a_{n}=a_0$ in which
case $A=S$.
\end{defn}

By considering a common refinement of two partitions $F_1$ and
$F_2$ of an arc $A\subset S$ such that $f(F_1)\cup f(F_2)\subset
T(S)$ and satisfying the conditions in Definition~\ref{partition},
it follows from Lemma~\ref{summ} that we get the same value for
$\var(f,A,S)$ whether we use the partition $F_1$ or the partition
$F_2$. Hence, $\var(f,A,S)$ is well-defined. If $A=S$ we denote
\index{variation@$\var(f,S)$}
$\var(f,S,S)$ simply by $\var(f,S)$\index{variation!of a simple closed curve}.

\subsection{Index and variation for finite
partitions}\label{indvar}

 What links Theorem~\ref{fpthm} with
variation is Theorem~\ref{I=V+1} below, first announced by Bell
 in the mid 1980's (see also Akis \cite{akis99}).  Our
proof is a modification of Bell's unpublished proof.  We first
need a variant of Proposition~\ref{straightjunction}. Let
$r:\complex\to T(S^1)$ be radial retraction: $r(z)=\frac{z}{|z|}$
when $|z|\geq 1$ and $r|_{T(S^1)}=id|_{T(S^1)}$.

\begin{lem} [Curve Straightening] \label{cs} Suppose $f:S^1\to\complex$ is a map with no
fixed points on $S^1$.  If $[a,b]\subset S^1$ is a proper subarc
with $f([a,b])\cap [a,b]=\0$, $f((a,b))\subset \complex\sm T(S^1)$
and $f(\{a,b\})\subset S^1$, then there exists a map
$\tf:S^1\to\complex$ such that $\tf|_{S^1\setminus (a,b)}=f|_{S^1\setminus (a,b)}$,
$\tf|_{[a,b]}:[a,b]\to (\complex\setminus T(S^1))\cup\{f(a),f(b)\}$ and $\tf|_{[a,b]}$ is
homotopic to $f|_{[a,b]}$ in $\{a,b\}\cup \complex\sm T(S)$ relative
to $\{a,b\}$, so that $r|_{\tf([a,b])}$ is   locally one-to-one.
Moreover, $\var(f,[a,b],S^1)=\var(\tf,[a,b],S^1)$.
\end{lem}

Note that if $\var(f,[a,b],S^1)=0$, then $r$ carries
$\tf([a,b])$ one-to-one onto the arc (or point) in $S^1\setminus (a,b)$ from
$f(a)$ to $f(b)$. If the $\var(f,[a,b],S^1)=m>0$, then
 $r\circ \tf$ wraps the arc $[a,b]$ counterclockwise about $S^1$ so that $\tf([a,b])$ meets each ray in $J_v$ $m$ times.
  A similar statement holds for negative variation. Note also
  that it is possible for index to be defined yet variation not
to be defined on a simple closed curve  $S$. For example, consider
the map $z\to 2z$ with $S$ the unit circle since there is no partition of $S$ satisfying the conditions in
Definition~\ref{vararc}.
\begin{thm}[Index = Variation + 1, Bell] \label{I=V+1} \index{index!Index=Variation+1 Theorem}
Suppose $g:S^1\to\complex$ is an
orientation preserving embedding onto a simple closed curve $S$
and $f:S\to\complex$ is a  fixed point free map.  If $F=\{a_0<a_1<\dots<a_n\}$
is a partition of $S$ and $A_i=[a_i,a_{i+1}]$ for $i=0,1,\dots,n$
with $a_{n+1}=a_0$ such that $f(F)\subset T(S)$ and $f(A_i)\cap
A_i=\0$ for each $i$, then
$$\ind(f,S)=\ind(f,g)=\sum_{i=0}^n\var(f,A_i,S)+1=\var(f,S)+1.$$
\end{thm}

\begin{proof}

By an appropriate conjugation of $f$ and $g$, we may assume
without loss of generality that $S=S^1$ and $g=id$.  Let $F$ and
$A_i=[a_i,a_{i+1}]$ be  as in the hypothesis. Consider the
collection of arcs
$$\mc{K}=\{K\subset S\mid \text{$K$ is the
closure of a component of $S\cap f^{-1}(f(S)\sm T(S))$}\}.$$
  For each
$K\in\mc{K}$, there is an $i$ such that $K\subset A_i$.  Since
$f(A_i)\cap A_i=\0$, it follows from the remark after
Definition~\ref{vararc} that
 $\var(f,A_i,S)=\sum_{K\subset A_{i}, K\in\mc{K}} \var(f,K,S)$.
 By the remark following Proposition~\ref{straightjunction}, we can compute $\var(f,K,S)$
 using one fixed junction for $A_i$. It is now clear that there are at most finitely many  $K\in\mc{K}$
 with $\var(f,K,S)\not=0$.
Moreover, the images of the endpoints of each $K$ lie on $S$.

Let $m$ be the cardinality of the set $\mc{K}_f=\{K \in \mc{K}
\mid
  \var(f,K,S)\not= 0\}$. By the above remarks, $m<\infty$ and  $\mc{K}_f$ is independent of
  the partition $F$.
  We prove the theorem by induction on $m$.

Suppose for a given $f$  we have $m=0$.  Observe that from the
definition of variation and the fact that the computation of
variation is independent of the choice of an appropriate
partition, it follows that,
$$\var(f,S)=\sum_{K\in\mc{K}} \var(f,K,S)=0.$$

We claim that there is a map $f_1:S\to\complex$ with
$f_1(S)\subset T(S)$ and a homotopy $H$ from $f|_S$ to $f_1$ such
that each level $H_t$ of the homotopy is fixed point free and
$\ind(f_1,id|_S)=1$.

To see the claim, first apply the Curve Straightening
Lemma~\ref{cs} to each $K\in\mc{K}$ (if there are infinitely many,
they form a null sequence) to obtain a fixed point free homotopy
of $f|_S$ to a map $\tf:S\to\complex$ such that
$r|_{\tf(K)}$ is
locally one-to-one on each $K\in\mc{K}$, where $r$ is radial
retraction  of $\complex$ to $T(S)$, and $\var(\tf,K,S)=0$ for each
$K\in\mc{K}$. Let $K$ be in $\mc{K}$ with endpoints $x,y$. Since
$\tf(K)\cap K=\0$ and $\var(\tf,K,S)=0$, $r|_{\tf(K)}$ is one-to-one, and $r\circ \tf(K)\cap K=0$.
Define $f_1|_K=r\circ \tf|_K$. Then $f_1|_K$ is fixed point free homotopic
 to $f|_K$ (with endpoints of $K$ fixed). Hence, if $K\in\mc{K}$ has
 endpoints $x$ and $y$, then
$f_1$ maps $K$ to the subarc of $S$ with endpoints $f(x)$ and $f(y)$ such that $K\cap f_1(K)=\0$.
Since $\mc{K}$ is a null family,
 we can do this for each $K\in\mc{K}$ and set $f_1|_{S^1\setminus \cup\mc{K}}=f|_{S^1\setminus \cup\mc{K}}$
  so that  we obtain the
desired $f_1:S\to\complex$ as the end map of a fixed point free
homotopy from $f$ to $f_1$.  Since $f_1$ carries $S$ into $T(S)$,
Corollary~\ref{mapinhull} implies $\ind(f_1,id|_S)=1$.

Since the homotopy $f\simeq f_1$ is fixed point free, it follows
from Theorem~\ref{fpfhomotopy2} that $\ind(f,id|_S)=1$. Hence, the
theorem holds if $m=0$ for any $f$ and any appropriate partition
$F$.

By way of contradiction suppose the collection $\mathcal{F}$ of all maps $f$
on $S^1$
 which satisfy the hypotheses of the theorem, but
not the conclusion is non-empty. By the above $0<|\mc{K}_f|<\infty$ for each.
Let $f\in\mathcal{F}$  be a counterexample for which
$m=|\mc{K}_f|$ is minimal.  By modifying $f$, we will show there
exists  $f_1\in\mathcal{F}$ with $|\mc{K}_{f_1}|<m$, a
contradiction.

Choose $K\in\mc{K}$ such that $\var(f,K,S)\not=0$.  Then
$K=[x,y]\subset A_i=[a_i,a_{i+1}]$ for some $i$.  By the Curve
Straightening Lemma~\ref{cs} and Theorem~\ref{fpfhomotopy2}, we
may suppose $r|_{f(K)}$ is locally one-to-one on $K$.  Define a new
map $f_1:S\to\complex$ by setting $f_1|_{\ol{S\sm K}}=f|_{\ol{S\sm
K}}$ and setting $f_1|_K$ equal to the linear map taking $[x,y]$
to the subarc $f(x)$ to $f(y)$ on $S$ missing $[x,y]$.
Figure~\ref{varfig} (left) shows an example of a (straightened)
$f$ restricted to $K$ and the corresponding $f_1$ restricted to $K$ for a case where $\var(f,K,S)=1$,
while Figure~\ref{varfig} (right) shows a case where
$\var(f,K,S)=-2$.

\begin{figure}

\includegraphics{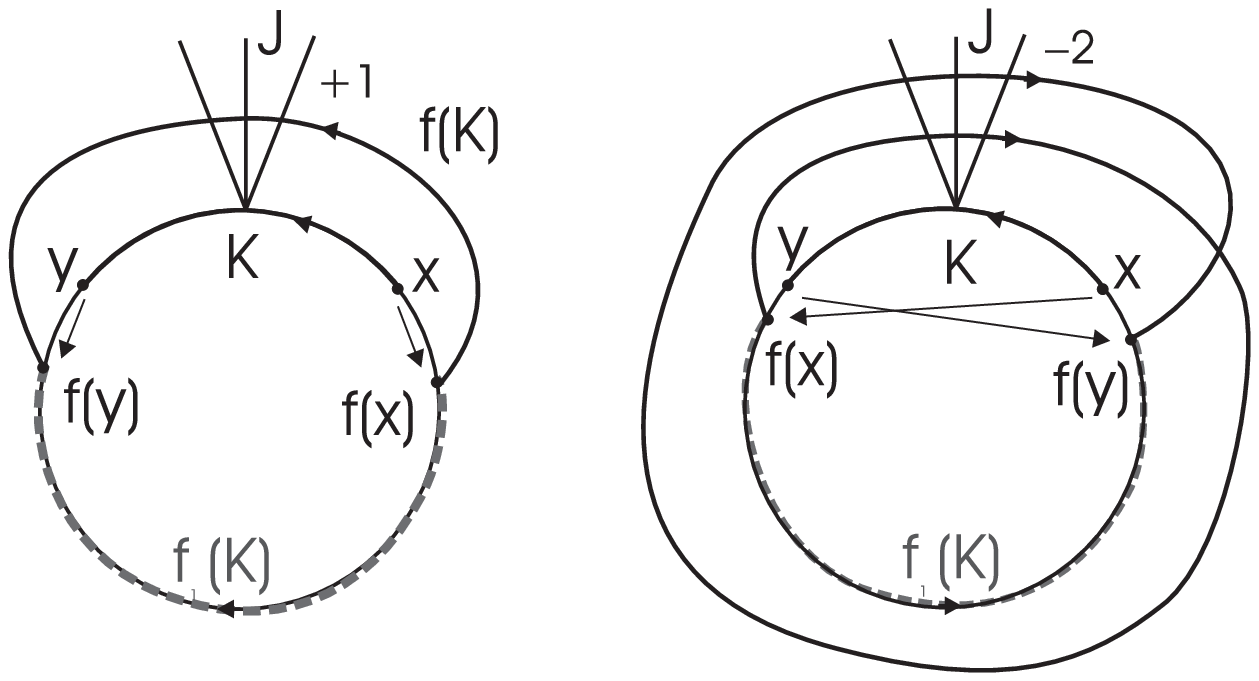}

\caption{Replacing $f:S\to\complex$ by $f_1:S\to\complex$ with one
less subarc of nonzero variation.} \label{varfig}

\end{figure}

Since on $\ol{S\sm K}$, $f$ and $f_1$ are the same map, we have
$$\var(f,S\sm K,S)=\var(f_1,S\sm K,S).$$   Likewise for the
fractional index, $$\ind(f,S\sm K)=\ind(f_1,S\sm K).$$ By
definition (refer to the observation we made in the case $m=0$),
$$\var(f,S)=\var(f,S\sm K,S)+\var(f,K,S)$$
$$\var(f_1,S)=\var(f_1,S\sm K,S)+\var(f_1,K,S)$$
and by Proposition~\ref{fracindex},
$$\ind(f,S)=\ind(f,S\sm K)+\ind(f,K)$$
$$\ind(f_1,S)=\ind(f_1,S\sm K)+\ind(f_1,K).$$ Consequently,
$$\var(f,S)-\var(f_1,S)=\var(f,K,S)-\var(f_1,K,S)$$ and
$$\ind(f,S)-\ind(f_1,S)=\ind(f,K)-\ind(f_1,K).$$

We will now show that the changes in index and variation, going
from $f$ to $f_1$ are the same (i.e., we will show that
$\var(f,K,S)-\var(f_1,K,S)=\ind(f,K)-\ind(f_1,K)$). We suppose
first that $\ind(f,K)=n+\alpha$ for some nonnegative $n\in\mathbb{N}$
and $0\le\alpha<1$. That is, the vector $f(z)-z$ turns through $n$
full revolutions counterclockwise and $\alpha$ part of a
revolution counterclockwise as $z$ goes from $x$ to $y$ counterclockwise along $S$.
(See Figure~\ref{varfig} (left) for a case $n=0$ and $\alpha$
about $0.8$.)  Then as $z$ goes from $x$ to $y$ counterclockwise along $S$, $f_1(z)$
goes along $S$ from $f(x)$ to $f(y)$ in the clockwise direction,
so $f_1(z)-z$ turns through  $-(1-\alpha)=\alpha-1$ part of a revolution.
Hence, $\ind(f_1,K)=\alpha-1$.  It is easy to see that
$\var(f,K,S)=n+1$ and $\var(f_1,K,S)=0$.  Consequently,
$$\var(f,K,S)-\var(f_1,K,S)=n+1-0=n+1$$ and
$$\ind(f,K)-\ind(f_1,K)=n+\alpha-(\alpha-1)=n+1.$$

In Figure~\ref{varfig} on the left we assumed that
$f(x)<x<y<f(y)$. The cases where $f(y)<x<y<f(x)$ and $f(x)=f(y)$
are treated similarly. In this case $f_1$ still wraps around in the positive direction,
but the computations are slightly different:
$\var(f,K)=1$, $\ind(f,K)=1+\al$, $\var(f_1,K)=0$ and $\ind(f_1,K)=\al$.

Thus when $n\geq 0$, in going from $f$ to $f_1$, the change in
variation and the change in index are the same.  However, in
obtaining $f_1$ we have removed one $K\in\mc{K}_f$, reducing the
minimal $m=|\mc{K}_f|$ for $f$ by one, producing a counterexample $f_1$
with $|\mc{K}_{f_1}|=m-1$, a contradiction.

The cases where $\ind(f,K)=n+\alpha$ for negative $n$ and
$0<\alpha<1$ are handled similarly, and illustrated for $n=-2$ and
$\alpha$ about $0.4$ in Figure~\ref{varfig} (right).

\end{proof}

\subsection{Locating arcs of negative variation}  The principal tool in
proving  Theorem~\ref{outchannel} (unique outchannel)  is the following
theorem first obtained by Bell (unpublished).  It provides a method for locating
arcs of negative variation on a curve of index zero.

\begin{thm} [Lollipop Lemma, Bell]  \label{lollipop} \index{Lollipop Lemma} Let $S\subset\complex$ be a
simple closed curve  and $f:T(S)\to\complex$ a fixed point free map.
Suppose $F=\{a_0<\dots<a_n<a_{n+1}<\dots<a_m\}$ is a partition of
$S$, $a_{m+1}=a_0$ and $A_i=[a_i,a_{i+1}]$   such that
$f(F)\subset T(S)$ and $f(A_i)\cap A_i=\0$ for $i=0,\dots, m$.
Suppose $I$ is an arc in $T(S)$ meeting $S$ only at its endpoints
$a_0$ and $a_{n+1}$. Let $J_{a_0}$ be a junction in $(\complex\sm
T(S))\cup\{a_0\}$ and suppose that $f(I)\cap (I\cup J_{a_0})=\0$.
Let $R=T([a_0,a_{n+1}]\cup I)$ and $L=T([a_{n+1},a_{m+1}]\cup I)$.
Then  one of the following holds:
\begin{enumerate} \item \label{first} If $f(a_{n+1})\in R$, then
$$\sum_{i\leq n}\var(f,A_i,S)+1=\ind(f,I\cup [a_0,a_{n+1}]) .$$

\item \label{second} If $f(a_{n+1})\in L$, then
$$\sum_{i>n}\var(f,A_i,S)+1=\ind(f,I\cup [a_{n+1},a_{m+1}]).$$

\end{enumerate}
\end{thm}
 (Note that in (\ref{first}) in effect we compute $\var(f,\bd R)$ but technically,
we have not defined $\var(f,A_i,\bd R)$  since the endpoints of
$A_i$ do not have to map inside $R$ but they do map into  $T(S)$.
Similarly in Case (\ref{second}).)
\begin{proof}

Without loss of generality, suppose $f(a_{n+1})\in L$.  Let
$C=[a_{n+1},a_{m+1}]\cup I$ (so $T(C)=L$). We want to construct a
map $f':C\to\complex$, fixed point free homotopic to $f|_C$, that
does not change variation on any arc $A_i$ in $C$ and has the
properties listed below.
\begin{enumerate}
\item \label{inL} $f'(a_i)\in L$ for all $n+1\le i\le m+1$. Hence
$\var(f',A_i,C)$ is defined for each $i>n$. \item \label{varf=g}
$\var(f',A_i,C)=\var(f,A_i,S)$ for all $n+1\le i\le m$. \item
\label{varI} $\var(f',I,C)=\var(f,I,S)=0$. \item \label{indf=g}
$\ind(f',C)=\ind(f,C)$.
\end{enumerate}
Having such a map, it then follows from Theorem~\ref{I=V+1}, that
$$\ind(f',C)=\sum_{i=n+1}^m \var(f',A_i,C)+\var(f',I,C)+1.$$

By Theorem~\ref{fpfhomotopy2} $\ind(f',C)=\ind(f,C)$.  By (\ref{varf=g})
and (\ref{varI}), $\sum_{i> n}
\var(f',A_i,C)+\var(f',I,C)=\sum_{i> n} \var(f,A_i,S)$ and  the
Theorem would follow.

\begin{figure}

\includegraphics{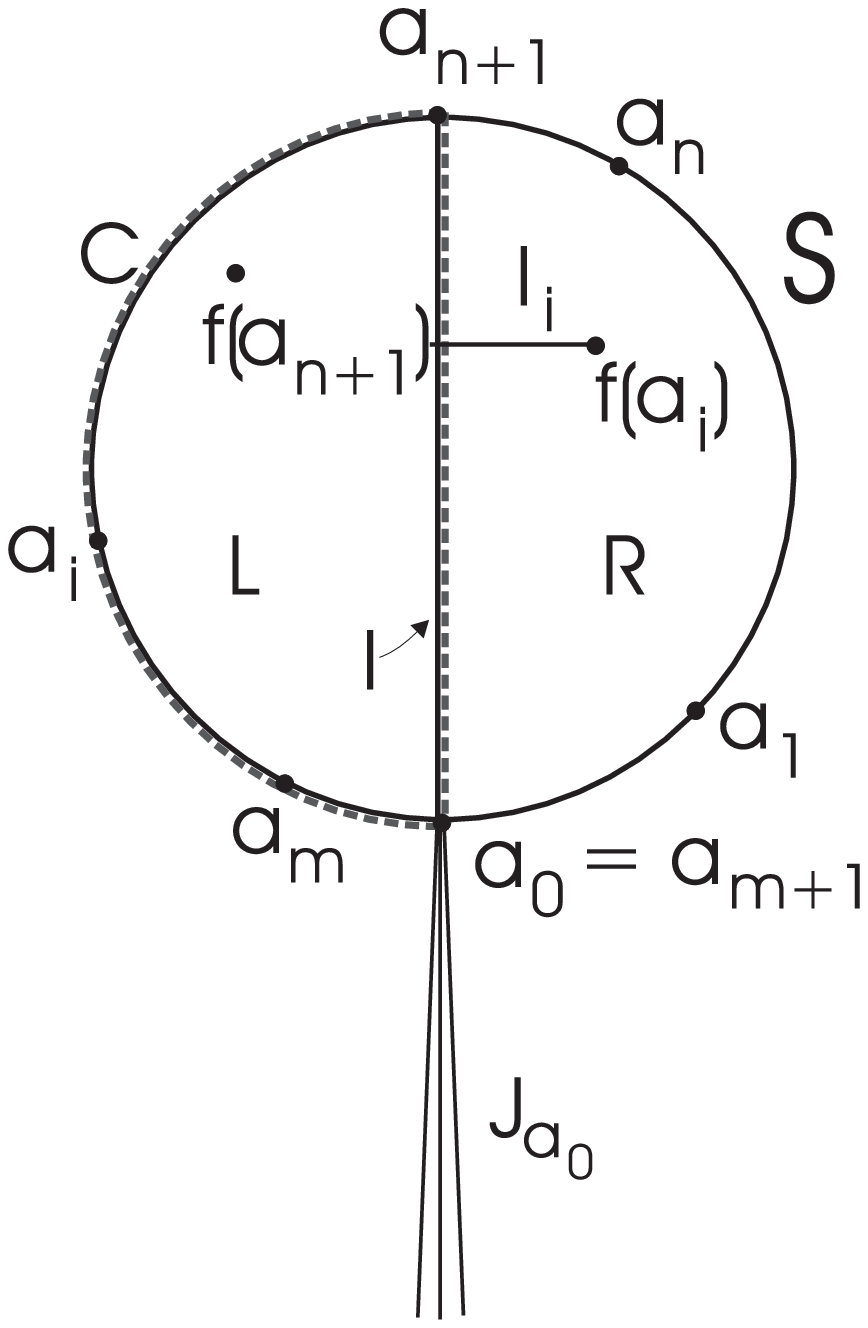}

\caption{Bell's Lollipop.} \label{lollypic}

\end{figure}

It remains to define the map $f':C\to\complex$ with the above
properties.  For each $i$ such that $n+1\le i \le m+1$, chose an
arc $I_i$ joining $f(a_i)$ to $L$ as follows:
\begin{enumerate}
\item[(a)] If $f(a_i)\in L$, let $I_i$ be the degenerate arc
$\{f(a_i)\}$. \item[(b)] If $f(a_i)\in R$ and $n+1<i<m+1$, let $I_i$
be an arc in $R\setminus\{a_0,a_{n+1}\}$ joining $f(a_i)$ to $I$.
\item[(c)] If $f(a_0)\in R$, let $I_0$ be an arc joining $f(a_0)$
to $L$ such that $I_0\cap(L\cup J_{a_0})\subset A_{n+1}\setminus
\{a_{n+1}\}$.
\end{enumerate}

Let $x_{n+1}=y_{n+1}=a_{n+1}$, $y_0=y_{m+1}\in
I\setminus\{a_0,a_{n+1}\}$ and $x_0=x_{m+1}\in
A_m\setminus\{a_m,a_{m+1}\}$. For $n+1<i<m+1$, let $x_i\in
A_{i-1}$ and $y_i\in A_i$ such that $y_{i-1}<x_i<a_i<y_i<x_{i+1}$.
For $n+1<i<m+1$ let $f'(a_i)$ be the endpoint of $I_i$ in $L$,
$f'(x_i)=f'(y_i)=f(a_i)$ and extend $f'$ continuously from
$[x_i,a_i]\cup [a_i,y_i]$ onto $I_i$ and define $f'$ from
$[y_i,x_{i+1}]\subset A_{i}$ onto $f(A_i)$ by
$f'|_{[y_i,x_{i+1}]}=f\circ h_i$, where $h_i:[y_i,x_{i+1}]\to A_i$
is a homeomorphism such that $h_i(y_i)=a_i$ and
$h_i(x_{i+1})=a_{i+1}$. Similarly, define $f'$ on
$[y_0,a_{n+1}]\subset I$ to $f(I)$ by $f|_{[y_0,a_{n+1}]}=f\circ
h_0$, where $h_0:[y_0,a_{n+1}]\to I$ is an onto homeomorphism such that
$h(a_{n+1})=a_{n+1}$ and extend  $f'$ from
$[x_{m+1},a_0]\subset A_m$ and $[a_o,y_0]\subset I$ onto $I_0$
such that $f'(x_{m+1})=f'(y_0)=f(a_0)$ and $f'(a_0)$ is the
endpoint of $I_0$ in $L$.

Note that $f'(A_i)\cap A_i=\0$ for $i=n+1,\dots,m$ and $f'(I)\cap
[I\cup J_{a_{0}}]=\0$. To compute the variation of $f'$ on each of
$A_m$ and $I$ we can use the junction $J_{a_{0}}$ Hence
$\var(f',I,C)=0$ and, by the definition of $f'$ on $A_m$,
$\var(f',A_m,C)=\var(f,A_m,S)$. For $i=n+1,\dots,m-1$ we can use
the same junction $J_{v_{i}}$ to compute $\var(f',A_i,C)$ as we
did to compute $\var(f,A_i,S)$. Since $I_i\cup I_{i+1}\subset
T(S)\sm A_i$ we have that $f'([a_i,y_i])\cup f'([x_{i+1},a_{i+1}])\subset
I_i\cup I_{i+1}$ misses that junction and, hence, make no
contribution to variation $\var(f',A_i,C)$. Since
$f'^{-1}(J_{v_{i}})\cap[y_i,x_{i+1}]$ is isomorphic to
$f^{-1}(J_{v_{i}})\cap A_i$, $\var(f',A_i,C)=\var(f,A_i,S)$ for
$i=n+1,\dots,m$.

To see that $f'$ is fixed point free homotopic to $f|_C$, note
that we can pull the image of $A_i$ back along the arcs $I_i$ and
$I_{i+1}$ in $R$ without fixing a point of $A_i$ at any level of
the homotopy.  Since $f'$ and $f|_C$ are fixed point free
homotopic and $f$ has no fixed points in $T(S)$, it follows from
Theorems~\ref{fpfhomotopy2} and \ref{fpthm}, that
$\ind(f',C)=\ind(f,C)=0$.
\end{proof}

Note that if $f$ is fixed point free on $T(S)$, then $\ind(f,S)=0$
and the next Corollary follows.
\begin{cor} \label{corlol}
Assume the hypotheses of Theorem~\ref{lollipop}. Suppose, in addition, $f$ is fixed point free on $T(S)$.
 Then if $f(a_{n+1})\in R$ there exists $i\leq n$ such that $\var(f,A_i,S)<0$.
 If $f(a_{n+1})\in L$ there exists $i> n$ such that $\var(f,A_i,S)<0$.
\end{cor}

\subsection{Crosscuts and bumping arcs}\label{infinitepartitions}\index{variation!for crosscuts}
For the remainder of Section 2,   our Standing Hypotheses are that
$f:\complex\to\complex$ takes continuum $X$ into $T(X)$ with no
fixed points in $T(X)$, and $X$ is minimal with respect to these
properties.

\begin{defn} [Bumping Simple Closed Curve] \label{bumping}\index{bumping!simple closed curve}
A simple closed curve $S$ in $\complex$ which has the property
that $S\cap X$ is nondegenerate and $T(X)\subset T(S)$ is said to
be a {\em bumping simple closed curve for $X$}. A subarc $A$ of a
bumping simple closed curve, whose endpoints lie in $X$, is said
to be a {\em bumping (sub)arc for $X$}.  \index{bumping!arc}Moreover, if $S'$ is any
bumping simple closed curve for $X$ which contains $A$, then $S'$
is said to \emph{complete} $A$. \index{completing a bumping arc}\end{defn}

A {\em crosscut} \index{crosscut} of $\tU=\rsphere\sm T(X)$ is an open arc $Q$
lying in $\tU$ such that $\ol Q$ is an arc with
endpoints $a\not=b\in T(X)$. In this case we will often write $Q=(a,b)$. (As seems to be traditional, we use
``crosscut of $T(X)$" interchangeably with ``crosscut of
$\tU$.") If $S$ is a bumping simple closed curve so that $X\cap S$ is
nondegenerate,
then each component of $S\sm X$ is a crosscut of $T(X)$. A similar
statement holds for a bumping arc $A$.
Given a non-separating continuum $T(X)$, let $A\subset\complex$ be a crosscut of $U^\infty=\rsphere\sm T(X)$.
Given a crosscut $A$ of $U^\infty$  denote by
$\Sh(A)$, the \emph{shadow of} \index{crosscut!shadow} \index{shadow}
\index{shadow@$\Sh(A)$}
 $A$,  the bounded component of $\complex\setminus [T(X)\cup A]$.

Since $f$ has no fixed points in $T(X)$ and $X$ is compact, we can
choose a bumping simple closed curve $S$ in a small neighborhood of $T(X)$
such that all crosscuts in $S\sm X$ are small, have positive distance to their image and so that $f$ has no fixed points in $T(S)$.
Thus, we obtain the following corollary to Theorem~\ref{fpthm}.

\begin{cor} \label{bumpingscc} There is a bumping simple closed curve $S$ for $X$
such that $f|_{T(S)}$ is fixed point free; hence, by \ref{fpthm},
$\ind(f,S)=0$. Moreover, any bumping simple closed curve $S'$ for $X$ such
that $S'\subset T(S)$ has $\ind(f,S')=0$. Furthermore, any
crosscut $Q$ of $T(X)$ for which $f$ has no fixed points in $T(X\cup
Q)$ can be completed to a bumping simple closed curve $S$ for $X$ for
which $\ind(f,S)=0$. \end{cor}

\begin{prop} \label{varcross} Suppose $A$ is a bumping subarc for $X$.
If $\var(f,A,S)$ is defined for some bumping simple closed curve
$S$ completing $A$, then for any bumping simple closed curve $S'$
completing $A$, $\var(f,A,S)=\var(f,A,S')$. \end{prop}

\begin{proof}  Since $\var(f,A,S)$ is defined, $A=\cup_{i=1}^n A_i$, where each $A_i$ is a bumping arc with $A_i\cap f(A_i)=\0$ and $|A_i\cap A_j|\le 1$ if $i\ne j$. By the remark following
Definition~\ref{partition}, it suffices to assume that $A\cap f(A)=\0$.  Let $S$ and $S'$ be two bumping simple closed curves completing $A$ for which
variation is defined.  Let $J_a$ and $J_{a'}$ be junctions whereby $\var(f,A.S)$ and
$\var(f,A,S')$ are respectively computed.  Suppose first that both junctions lie (except
for $\{a,a'\}$) in $\complex\sm (T(S)\cup T(S'))$.  By the Junction Straightening
Proposition~\ref{straightjunction}, either junction can be used to compute either
variation on $A$, so the result follows. Otherwise, at least one junction is not in
$\complex\sm (T(S)\cup T(S'))$.  But both junctions are in $\complex\sm T(X\cup A)$.
Hence, we can find another bumping simple closed curve $S''$ such that $S''$ completes $A$,
 and
both junctions lie in $(\complex\sm T(S''))\cup\{a,a'\}$. Then by the
Propositions~\ref{straightjunction} and the definition of variation,
$\var(f,A,S)=\var(f,A,S'')=\var(f,A,S')$.
\end{proof}

It follows from Proposition~\ref{varcross} that  variation on a
crosscut $Q$, with $\ol{Q}\cap f(\ol{Q})=\0$, of $T(X)$ is independent of the bumping simple closed curve $S$
for $T(X)$ of which $Q$ is a subarc and is such that $\var(f,S)$ is defined. Hence, given a bumping arc $A$
of $X$,  we can denote
$\var(f,A,S)$ simply by $\var(f,A)$ \index{varfA@$\var(f,A)$}
when $X$ is understood.

The following proposition follows from Corollary~\ref{bumpingscc}, Proposition~\ref{varcross} and
Theorem~\ref{I=V+1}.

\begin{prop} \label{crossunder}  Suppose $Q$ is a crosscut of $T(X)$ such that
$f$ is fixed point free on $T(X\cup Q)$ and $f(\ol{Q})\cap \ol{Q}=\0$. Suppose $Q$ is replaced by
a bumping subarc $A$ with the same endpoints such that $Q\cup T(X)$ separates $A\sm X$ from
$\infty$ and each component $Q_i$ of $A\sm X$ is a crosscut such that $f(\cl{Q_i})\cap
\cl{Q_i}=\0$.  Then
$$\var(f,Q,X)=\sum_{i}
\var(f,Q_i,X)=\var(f,A,X).$$
\end{prop}

\subsection{Index and Variation for Carath\'eodory Loops}

We extend the definitions of index and variation
 to {\em Carath\'eodory loops}.

\begin{defn}[Carath\'eodory Loop]\index{Carath\'eodory Loop} Let $g:S^1\to\complex$ such
that $g$ is continuous and has a continuous
 extension $\psi:\ol{\rsphere\sm
T(S^1)}\to\ol{\rsphere\sm T(g(S^1))}$ such that
$\psi|_{\complex\setminus T(S^{1})}$ is an orientation preserving
homeomorphism from $\complex\setminus T(S^{1})$ onto
$\complex\setminus T(g(S^1))$. We call $g$ (and loosely,
$S=g(S^1)$), a {\em Carath\'eodory loop}.
\end{defn}

In particular, if $g:S^1\to g(S^1)=S$ is a continuous extension of a Riemann map $\psi:\disk^\infty\to\rsphere\sm T(g(S^1))$, then $g$ is
 a Carath\'eodory loop, where $\disk^\infty=\{z\in\rsphere\mid |z|>1\}$ \index{deltai@$\disk^\infty$}
 is the ``unit disk'' about $\infty$.

Let $g:S^1\to\complex$ be a Carath\'eodory loop and let
$f:g(S^1)\to\complex$ be a fixed point free map.  In order to define variation of $f$
on $g(S^1)$, we do the partitioning in
$S^1$ and transport it to the Carath\'eodory loop $S=g(S^1)$.
An {\em allowable} partition \index{allowable partition} of $S^1$ is a set
$\{a_0<a_1<\dots<a_n\}$ in $S^1$ ordered counterclockwise, where
$a_0=a_n$ and $A_i$ denotes the counterclockwise interval
$[a_{i},a_{i+1}]$, such that for each $i$, $f(g(a_i))\in
T(g(S^1))$ and $f(g(A_i))\cap g(A_i)=\0$. Variation
$\var(f,A_i,g(S^1))=\var(f,A_i)$ on
each path $g(A_i)$ is then defined exactly as in
Definition~\ref{vararc}, except that the junction (see
Definition~\ref{junction}) is chosen so that the vertex
$v\in g(A_i)$ and $J_v\cap T(g(S^1))\subset \{v\}$, and the
crossings of the junction $J_v$  by $f(g(A_i))$ are counted (see
Definition~\ref{vararc}). Variation on the whole loop, or an
allowable subarc thereof, is defined just as in
Definition~\ref{partition}, by adding the variations on the
partition elements. At this point in the development, variation is
defined only relative to the given allowable partition $F$ of
$S^1$ and the parameterization $g$ of $S$:
$\var(f,F,g(S^1))$.

Index on a Carath\'eodory loop $S$ is defined exactly as in
Section~\ref{defindex} with $S=g(S^1)$ providing the
parameterization of $S$. Likewise, the definition of fractional
index and Proposition~\ref{fracindex} apply to Carath\'eodory
loops.

Theorems~\ref{fpfhomotopy}, \ref{fpfhomotopy2},
Corollary~\ref{mapinhull}, and Theorem~\ref{fpthm} (if $f$ is also defined on $T(S)$)
apply to Carath\'eodory loops.  It follows that index on a Carath\'eodory
loop $S$ is independent of the choice of parameterization $g$.
The Carath\'eodory loop $S$ is approximated, under small homotopies, by simple closed curves
$S_i$. Allowable partitions of $S$ can be made to correspond to allowable partitions
of $S_i$ under small homotopies. Since variation and index are invariant under suitable
homotopies (see the comments after Proposition~\ref{straightjunction}) we have the following theorem.

\begin{thm} \label{CaraVI} \index{index!I=V+1 for Carath\'eodory Loops}
Suppose $S=g(S^1)$ is a parameterized
Carath\'eodory loop in $\complex$ and $f:S\to\complex$ is a fixed point
free map. Suppose variation of $f$ on $S^1=A_0\cup\dots\cup A_n$ with
respect to $g$ is defined for some partition $A_0\cup\dots\cup
A_n$ of $S^1$. Then
$$\ind(f,g)=\sum_{i=0}^n\var(f,A_i,g(S^1))+1.$$
\end{thm}

\subsection{Prime Ends}
Prime ends provide a way of studying the approaches to the boundary of a
simply-connected plane domain with non-degenerate boundary.  See \cite{colllohw66} or
\cite{miln00} for an analytic summary of the topic and \cite{urseyoun51} for a more
topological approach. We will be interested in the prime ends of $\tU=\rsphere\sm
T(X)$.  Recall that  $\disk^\infty=\{z\in\rsphere\mid |z|>1\}$ is the ``unit disk
about $\infty$."  The Riemann Mapping Theorem guarantees the existence of a conformal
map $\phi:\disk^\infty\to \tU$ taking $\infty\to\infty$, unique up to the argument
of the derivative at $\infty$. Fix such a map $\phi$. We identify $S^1= \bd\disk^\infty$
with $\real/\zed$ and identify points $e^{2\pi i t}$ in $\bd\disk^\infty$ by their
argument $t \pmod{1}$. Crosscut and shadow were defined in Section~\ref{infinitepartitions}.

\begin{defn}[Prime End]  A {\em chain of crosscuts}\index{chain of crosscuts} \index{prime end}is a
sequence $\{Q_i\}_{i=1}^\infty$ of crosscuts of $\tU$ such
that for $i\not= j$, $Q_i\cap Q_j=\0$, $\mrm{diam}(Q_i)\to 0$, and
for all $j>i$, $Q_i$ separates $Q_j$ from $\infty$ in $\tU$. Hence, for all $j>i$,
$Q_j\subset \Sh(Q_i)$.
Two chains of crosscuts are said to be {\em equivalent} \index{chain of crosscuts!equivalent} iff it is
possible to form a sequence of crosscuts by selecting alternately
a crosscut from each chain so that the resulting sequence of
crosscuts is again a chain. A {\em prime end} $\mc{E}$ is an
equivalence class of chains of crosscuts.\index{Et@$\mc{E}_t$}\end{defn}

If $\{Q_i\}$ and $\{Q'_i\}$  are equivalent  chains of crosscuts of $\tU$, it can be
shown that $\{\phi^{-1}(Q_i)\}$ and $\{\phi^{-1}(Q'_i)\}$ are equivalent  chains of crosscuts of
$\disk^\infty$  each of which converges to the same unique point $e^{2\pi it}\in S^1=
\bd\disk^\infty$, $t\in[0,1)$, independent of the representative chain. Hence, we denote by
 $\mc{E}_t$ the prime end of $\tU$ defined by $\{Q_i\}$.

\begin{defn}[Impression and Principal Continuum] \index{prime end!impression}\index{prime end!principal continuum}
\index{impression} \index{principal continuum}
Let $\mc{E}_t$ be a prime end of $\tU$ with defining chain of crosscuts
$\{Q_i\}$.
The set \index{imet@$\im(\mc{E}_t)$}
$$\im(\mc{E}_t)=\bigcap_{i=1}^\infty \cl{\Sh(Q_i)}$$ is a subcontinuum
of $\bd \tU$ called the {\em impression} of $\mc{E}_t$.  The
set \index{pret@$\pr(\mc{E}_t)$} $$\pr(\mc{E}_t)=\{z\in\bd \tU\mid \text{for some chain
$\{Q'_i\}$ defining $\mc{E}_t$, $Q'_i\to z$}\}$$ is a continuum
called the {\em principal continuum} of $\mc{E}_t$.
\end{defn}

For a prime end $\mc{E}_t$, $\pr(\mc{E}_t)\subset \im(\mc{E}_t)$,
possibly properly.    We will be interested in the existence of
prime ends $\mc{E}_t$ for which $\pr(\mc{E}_t)=\im(\mc{E}_t)=\bd
\tU$.

\begin{defn} [External Rays] \index{external ray}  \index{Rt@$R_t$} Let $t\in [0,1)$ and define
$$R_t=\{z\in\complex\mid z=\phi(re^{2\pi it}),1<r<\infty\}.$$ We
call $R_t$ the {\em external ray (with argument $t$)}. If $x\in R_t$ then the
$(X,x)$-\emph{end of} $R_t$ is the  bounded component $K_x$ of
\index{external ray!end of}
$R_t\sm\{x\}$.\end{defn}

The external rays $R_t$ foliate $\tU$.

\begin{defn}[Essential crossing]\label{essential}\index{external ray!essential crossing}
An external ray $R_t$ is said to \emph{cross} a crosscut $Q$
\emph{essentially} if and only if there exists $x\in R_t$ such that the $(T(X),x)$-end of $R_t$ is
contained in the bounded complementary domain of $T(X)\cup Q$. In this case we will also say that $Q$ crosses $R_t$ essentially.
\end{defn}

The results listed below are known.
\begin{prop}[\mbox{\cite{colllohw66}}] \label{trans}
Let $\mc{E}_t$ be a prime end of $\tU$. Then
$\pr(\mc{E}_t)=\cl{R_t}\sm R_t$.  Moreover, for each $1<r<\infty$
there is a  crosscut $Q_r$ of $\tU$ with  $\{\phi(re^{2\pi it})\}=R_t\cap Q_r$ and
$\dm(Q_r)\to 0$ as $r\to 1$ and such that $R_t$ crosses $Q_r$
essentially.
\end{prop}

\begin{defn}[Landing Points and Accessible Points]  \index{accessible point}\index{external ray!landing point}If
$\pr(\mc{E}_t)=\{x\}$, then we say $R_t$ {\em lands} on $x\in
T(X)$ and $x$ is the {\em landing point} \index{landing point} of $R_t$.  A point $x\in
\bd T(X)$ is said to be {\em accessible} (from $\tU$) iff
there is an arc in $\tU\cup\{x\}$ with
$x$ as one of its endpoints. \end{defn}

\begin{prop} A point $x\in \bd T(X)$ is accessible iff $x$ is the
landing point of some external ray $R_t$. \end{prop}

\begin{defn}[Channels]\index{channel}\index{prime end!channel}
A prime end $\mc{E}_t$ of $\tU$ for which $\pr(\mc{E}_t)$ is
nondegenerate is said to be a {\em channel in} $\bd \tU$ (or
in $T(X)$). If moreover $\pr(\mc{E}_t)=\bd \tU=\bd T(X)$, we
say $\mc{E}_t$ is a {\em dense} channel. \index{channel!dense} A crosscut $Q$ of $\tU$
 is said to {\em cross} the channel $\mc{E}_t$
iff  $R_t$ crosses $Q$ essentially.
\end{defn}

When $X$ is locally connected, there are no channels, as the
following classical theorem proves.  In this case, every prime end
has degenerate principal set and degenerate impression.

\begin{thm} [Carath\'eodory] \label{carath} $X$ is locally connected iff the
Riemann map $\phi:\disk^\infty\to \tU=\rsphere\sm T(X)$
taking $\infty\to\infty$ extends continuously to
$S^1=\bd\disk^\infty$.
\end{thm}

\section{Kulkarni-Pinkall Partitions}\label{secKP} Throughout this section
let $K$ be a compact subset of the plane whose complement
$U=\complex\sm K$ is connected. In the interest of completeness we define the Kulkarni-Pinkall
partition of $U$ and prove the basic properties of this partition that are essential for
our work in Section~\ref{sechyp}. Kularni-Pinkall \cite{kulkpink94} worked in closed $n$-manifolds.
We will follow their approach and adapt it to our situation in the plane.

We think of $K$ as a closed subset of the Riemann sphere $\sphere$, with the spherical metric
and set $\tU=\sphere\sm K= U\cup\{\infty\}$. Let $\B^\infty$ \index{ball@$\B^\infty$} be the family of closed, round  balls $B$
in $\sphere$ such that $Int(B)\subset\tU$ and $|\partial B\cap K|\ge 2$. Then $\B^\infty$ is
in one-to-one correspondence with the family $\B$ \index{ball@$\B$} of closed subsets $B$ of $\complex$ which are the closure
of a complementary component of a straight line or a round circle in $\complex$ such that
$Int(B)\subset U$ and $|\partial B\cap K|\ge 2$.\index{maximal ball}

\begin{prop}\label{lense} If $B_1$ and $B_2$ are two closed round balls in $\complex$ such that
$B_1\cap B_2\ne \0$ but does not contain a diameter of either $B_1$ or $B_2$,
then $B_1\cap B_2$ is contained in a ball of diameter strictly less than the
diameters of both  $B_1$ and $B_2$.
\end{prop}

\begin{proof} Let $\partial B_1\cap \partial B_2=\{s_1,s_2\}$. Then the closed ball
with center $(s_1+s_2)/2$ and radius
$|s_1-s_2|/2$ contains  $B_1\cap B_2$.
\end{proof}

If $B$ is the closed ball of minimum diameter
that contains $K$, then we say that $B$ is the \textit{smallest ball} \index{smallest ball}
containing $K$.
It is unique by Proposition~\ref{lense}. It
exists, since any sequence of balls of decreasing diameters that
contain $K$ has a convergent subsequence. \medbreak

We denote the \emph{Euclidean convex hull of  $K$} \index{convex hull!Euclidean} by
$\ECH(K)$ \index{conve@$\ECH(K)$}.
 It is the intersection of all closed half-planes (a closed half-plane
is the closure of a component of the complement of a straight line) which contain
$K$. Hence $p\in\ECH(K)$ if $p$ cannot be separated from $K$ by a
straight line.

Given a closed ball $B\in\B^\infty$,  $\Int(B)$ is conformally  equivalent to
 the unit disk in $\complex$.
Hence its interior can be naturally equipped with the hyperbolic metric. Geodesics $\rg$
\index{g@$\rg$}  in
this metric are intersections of $\Int(B)$ with round circles $C\subset\rsphere$  which perpendicularly
cross the boundary $\partial B$. For every hyperbolic geodesic $\rg$, \index{hyperbolic!geodesic}
$B\setminus \ol{\rg}$ has exactly two components.
We call the closure of such components \emph{hyperbolic half-planes of
$B$.}\index{hyperbolic!halfplane} Given $B\in\B^\infty$, the \emph{hyperbolic convex hull of $K$ in $B$}\index{convex hull!hyperbolic}
 is the intersection of
 all (closed) hyperbolic
half-planes of $B$ which contain $K\cap B$ and we denote it by $\HCH(B\cap K)$. \index{convea@$\HCH(B\cap K)$}

\begin{lem}\label{nonsep}
Suppose that $B\in\B$ is the smallest ball containing $K\subset\complex$ and let $c\in B$
be its center. Then $c\in\HCH (K\cap\partial B)$.
\end{lem}
\begin{proof}
By contradiction. Suppose that there exists a circle that separates
the center $c$ from $K\cap\partial B$ and crosses $\partial B$
perpendicularly. Then there exists a line $\ell$ through $c$ such
that a half-plane bounded by $\ell$ contains $K\cap \partial B$ in
its interior. Let $B'=B+v$ be a translation of $B$ by a vector $v$
that is orthogonal to $\ell$ and directed into this halfplane. If $v$
is sufficiently small, then $B'$  contains $K$ in its interior. Hence, it can
be shrunk to a strictly smaller ball that also contains $K$, contradicting that $B$ has smallest
diameter.
\end{proof}

\begin{figure}
\includegraphics{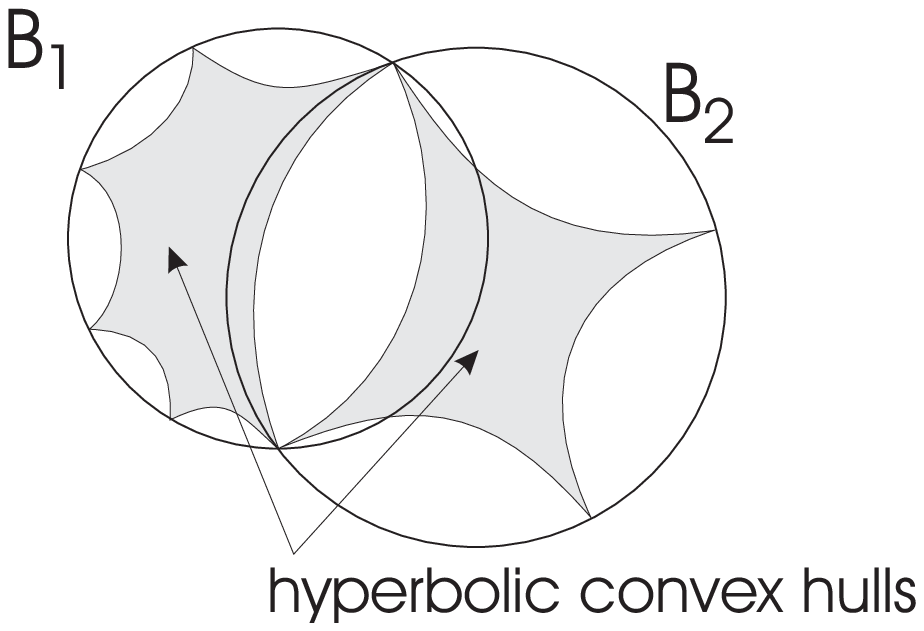}
\caption{Maximal balls have disjoint hulls.} \label{hulls}
\end{figure}
 \begin{lem}\label{disjointCH}
Suppose that $B_1,B_2\in\B^\infty$ with $B_1\ne B_2$.
Then
$$\HCH(B_1\cap\partial U)\cap\HCH(B_2\cap \partial
U)\subset\partial U.$$ In particular,
$\HCH(B_1\cap\partial U)\cap\HCH(B_2\cap \partial
U)$  contains at
most two points.
 \end{lem}
\begin{proof}
 A picture easily explains this, see Figure~\ref{hulls}. Note that $\partial U\cap [B_1\cup B_2]
 \subset \partial (B_1\cup B_1)$. Therefore $B_1\cap\partial U$
 and $B_2\cap\partial U$ share at most two points. The open hyperbolic chords
 between these points in the respective balls are disjoint.
\end{proof}

It follows that any point in $U^\infty$ can be contained in at most one
hyperbolic convex hull. In the next lemma we see that each point of $\tU$ is
indeed contained in $\HCH(B\cap K)$ for some $B\in\B^\infty$. So
$\{\tU\cap \HCH(B\cap K) \mid B\in\B^\infty\}$ is a
partition of $\tU$.

Since hyperbolic convex hulls are preserved by M\"{o}bius
transformations, they are more easy to manipulate than the Euclidean
convex hulls used by Bell  (which are preserved only by M\"{o}bius transformations
that fix $\infty$). This is illustrated by the proof of the following
lemma.

\begin{lem}[Kulkarni-Pinkall inversion lemma]\label{kplem} \index{Kulkarni-Pinkall!Lemma} For any
$p\in \sphere\setminus K$ there exists  $B\in\B^\infty$
such that $p\in\HCH(B\cap K)$.
\end{lem}

\begin{proof} We prove first that there exists $B^*\in\B^\infty$ such that
no line or circle which crosses $\partial B^*$ perpendicularly separates
$K\cap \partial B^*$ from $\infty$.

Let $B'$ be the smallest round ball which contains $K$ and let $B=\ol{\complex\sm B'}$.
Then $B^*=B\cup \{\infty\}\in \B^\infty$.
If $L$ is a  circle which crosses $\partial B^*=\partial B'$
perpendicularly and separates $K\cap \partial B'$ from $\infty$, then it also
separates $K\cap\partial B'$ from the center $c'$  of $B'$,
contrary to Lemma~\ref{nonsep}. [To see this note that if a hyperbolic geodesic $\rg$
of $B^*$ separates $K\cap B^*$ from $\infty$, then $\rg$ is contained in a round circle
$C$ and $C\cap B'$ separates $c'$ from $B'\cap K$, a contradiction.]  Hence,
$\infty \in\HCH(B^*\cap K)$.

Now let $p\in \complex^\infty\sm K$. Let $M:\sphere\to\sphere$ be a
M\"{o}bius transformation such that $M(p)=\infty$. By the above argument
there exists a ball  $B^*\in\B^\infty$ such that $\infty\in \HCH(B^*\cap M(K))$.
Then $B=M^{-1}(B^*)\in\B^\infty$ and, since $M$ preserves perpendicular circles, $p\in\HCH(B\cap K)$
as desired.
\end{proof}

From Lemmas~\ref{disjointCH} and \ref{kplem}, we obtain the following Theorem
which is a special case of a Theorem of  Kulkarni and
Pinkall \cite{kulkpink94}.

\begin{thm}\label{KPthrm}\index{Kulkarni-Pinkall!Partition}
 Suppose that $K\subset\complex$ is a nondegenerate compact set such that its complement
$\tU$ in the Riemann sphere is non-empty and connected. Then
$\tU$ is partitioned by the family
$$\KPP=\left\{\tU\cap \HCH(B\cap K)\colon B\in\B^\infty\right\}.$$
\end{thm}
Theorem~\ref{KPthrm} is the linchpin of the theory of geometric crosscuts.
An analogue of it  was known to Harold Bell and used by him implicitly since the
early 1970's. Bell considered non-separating plane continua $K$ and he
used the equivalent notion of Euclidean convex hull of the sets
$B\cap\partial U$ for all maximal balls $B\in\B$ (see
 the comment following Theorem~\ref{Hypmain}).

 \smallskip
 We will denote the partition $\left\{\tU\cap \HCH(B\cap K)\colon B\in\B^\infty\right\}$
of $\tU$ by $\KPP$. \index{KPP@$\KPP$}
Let $B\in \B^\infty$. If $B\cap\partial U^\infty$ consists of two points $a$ and $b$, then its (hyperbolic) hull \index{hull!hyperbolic}
is an open circular segment $\rg$ with endpoints $a$ and $b$ and perpendicular to $\partial B$.
 We will call the crosscut
$\rg$  a \textit{$\kp$
crosscut} or simply a \textit{$\kp$ chord}\index{KPchord@$\kp$ chord}. If $B\cap\partial U$ contains
three or more points, then we say that the hull  $\HCH(B\cap\partial
U)$ is a \textit{gap}.  \index{gap} A gap has nonempty interior. Its boundary in
$\Int(B)$ is a union of open circular segments (with endpoints in $K$), which we also call \emph{$\kp$
crosscuts or $\kp$ chords}. We denote by $\KP$ \index{KP@$\KP$}  the collection of all open  chords  obtained
as above using all  $B\in\B^\infty$.
\medbreak

The following example may serve to illustrate
Theorem~\ref{KPthrm}.\\
{\bf Example}.  Let $K$ be the unit square $\{x+yi\colon
-1\leq x,y\leq 1\}$.
There are five obvious members of $\B$. These are the sets
 $$\im
z\geq 1,\ \im z\leq -1,\ \re z\geq 1,\ \re z\leq -1,\ |z|\geq \sqrt
2,$$ four of which are half-planes. These are the only members of $\B$
whose hyperbolic convex hulls
have non-empty interiors. However,
for this example the family $\B$ defined
in the introduction of Section~\ref{secKP} is infinite.
 The hyperbolic hull of the half-plane $\im z\ge 1$ is the semi-disk
$\{z \mid  |z-i|\le 1, \ \im z >1\}$. The hyperbolic hulls of the other three half-planes given above
are also semi-disks. The hyperbolic hull of $|z|\ge \sqrt{2}$ is the unbounded region whose boundary
consists of the four semi-circles lying (except for their endpoints) outside $K$
and contained in the circles of radius $\sqrt{2}$ and having centers at $-2,2,-2i$
and $2i$, respectively. These hulls do not cover $U$ as there are spaces
between the hulls of the half-planes and the hull of $|z|\geq \sqrt
2$.

If $C$ is a circle that circumscribes $K$ and contains exactly two of its
vertices, such as $1\pm i$, then the exterior ball $B$ bounded by
$C$ is maximal. Now $\HCH(B\cap K)$ is a single chord and the union
of all such chords foliates the remaining spaces in $\complex\sm K$.

\begin{lem}\label{contKP}
If  $\rg_i$
is a sequence of $\kp$ chords with endpoints $a_i$ and $b_i$, and $\lim a_i=a\ne b=\lim b_i$,
then $\ol{\rg_i}$ is  convergent  and $\lim \ol{\rg_i}=C$, where $\rg=C\sm \{a,b\}\in\KP$ is also a
$\kp$ chord.
\end{lem}

\begin{proof}
For each $i$ let $B_i\in\B^\infty$ such that $\rg_i\subset\HCH(B_i\cap K)$.
Then $B_i$ converges to some $B\in\B^\infty$ and $\ol{\rg_i}$ converges to a closed
circular arc $C$ in $B$ with endpoints $a$ and $b$, and $C$ is perpendicular to $\partial B$.
Hence $\rg=C\sm\{a,b\}\subset \HCH(B\cap K)$. So $\rg\in\KP$.
\end{proof}

By Lemma~\ref{contKP}, the family $\KP$ of  chords has continuity
properties similar to a foliation.

\begin{lem}\label{lenseunion}
For $a,b\in K\cap\partial U^\infty$, define $C(a,b)$\index{Cab@$C(a,b)$} as the union of all $\kp$ chords  with endpoints
$a$ and $b$. Then if $C(a,b)\ne\0$, $C(a,b)$ is either  a single chord, or
$C(a,b)\cup\{a,b\}$ is a closed disk whose boundary consists of two $\kp$ chords contained in
$C(a,b)$ together with   $\{a,b\}$.
\end{lem}

\begin{proof}
Suppose $\rg$ and $\rh$ are two distinct $\kp$ chords between $a$ and $b$.
Then $S=\rg\cup\rh\cup\{a,b\}$ is a simple closed curve.
Choose a point $z$ in the   complementary domain $V$ of $S$ contained in $\tU$.
 Since the hyperbolic hulls  partition
$\tU$, there exists $B\in \B^\infty$ such that $z\in \HCH(B\cap K)$
and $\HCH(B\cap K)$ can only intersect $S\cap K$ in $\{a,b\}$. So
$\HCH(B\cap K)\cap K=\{a,b\}$ and it follows that
$V$ is contained in $C(a,b)$.

The rest of the Lemma follows from \ref{contKP}.
\end{proof}

\section{Hyperbolic foliation of simply connected domains}\label{sechyp}
In this section we will apply the results from Section~\ref{secKP} to the case that $K$
is a non-separating plane continuum (or, equivalently, that $\tU=\sphere\sm K$ is
simply connected). The results in this section are essential to \cite{overtymc07} but are
not used in this paper. The reader who is only interested in the fixed point question can skip this
section.

Let $\Disk$ \index{da@$\Disk$} be the open unit disk in the plane.
In this section we let $\phi:\Disk\to \sphere\sm K=\tU$ be a Riemann map onto $\tU$. We endow $\mathbb
D$ with the hyperbolic metric, which is carried to $\tU$ by the Riemann map.
We use $\phi$ and  the Kulkarni-Pinkall hulls to induce a
closed collection  $\Gamma$ of chords in $\mathbb D$ that is a hyperbolic geodesic
lamination in $\Disk$ (see \cite{thur85}).

Let $\rg\in\KP$ be a chord with endpoints $a$ and $b$. Then $a$ and $b$ are accessible points
in $K$ and $\ol{\phi^{-1}(\rg)}$ is an arc in $\Disk$ with endpoints $z,w\in\partial\Disk$.
 Let $\rG$ \index{G@$\rG$}  be the hyperbolic geodesic \index{geodesic!hyperbolic}
in $\Disk$ joining $z$ and $w$. Let $\Gamma$  \index{Gamma@$\Gamma$} be the collection of all $\rG$ such that $\rg\in\KP$.
 We will prove that  $\Gamma$  inherits  the properties of the family  $\KP$ as described
 in Theorem~\ref{KPthrm} and Lemma~\ref{contKP}
(see Lemma~\ref{compactfol}, Theorem~\ref{Hypmain} and the remark following ~\ref{Hypmain}).

 Since members of $\KP$  do not intersect (though their closures are arcs which may have common
 endpoints)  the same is true for distinct members of $\Gamma$. We will refer to the members of $\Gamma$
(and their images under $\phi$) as \emph{hyperbolic chords}
or \emph{hyperbolic geodesics} \index{hyperbolic chord}. Given $\rg\in\KP$ we denote
the corresponding element of $\Gamma$ by $\rG$ and its image  $\phi(\rG)$ in $\tU$ by $\fg$.
\index{g@$\fg$}
Note that $\Gamma$ is a lamination of $\mathbb{D}$ in the sense of Thurston\cite{thur85}.
By a \emph{gap} of $\Gamma$ (or of $\phi(\Gamma)$), we mean the closure of a component
of $\mathbb{D}\sm\bigcup \Gamma$ in $\mathbb{D}$ (or its image under $\phi$ in $\tU$, respectively).

\begin{lem}[J{\o}rgensen~{\cite[p.91]{pomm92}}]\label{jorg}\index{J{\o}rgensen Lemma}
Let $B$ be a closed round ball such that its interior is in $\tU$. Let
$\gamma\subset\mathbb D$ be a hyperbolic geodesic.
Then $\phi(\gamma)\cap B$ is connected.
In particular, if $R_t$ is an external ray in $U^\infty$ and $B\in\B^\infty$, then
$R_t\cap B$ is connected.
\end{lem}

If $a,b\in\partial\tU$,  recall that $C(a,b)$ is the union of all $\kp$ chords  with endpoints
$a$ and $b$. From the viewpoint of prime ends, all chords in $C(a,b)$ are the
same. That is why all the chords in $C(a,b)$ are
replaced  by a single hyperbolic
chord  $\fg\in\phi(\Gamma)$.
The following lemma follows.

\begin{lem}\label{sameb}Suppose $\rg\in\kp$ and $\rg\subset \HCH(B\cap \partial\tU)$ joins the points
$a,b\in\partial\tU$ for some $B\in B^\infty$. We may assume that the Riemann map
$\phi:\Disk\to\tU$ is extended over all points $x\in S^1$ so that $\phi(x)$ is an accessible point
of $\tU$. Let $\phi^{-1}(a)=\ta$, $\phi^{-1}(b)=\tb$ and $\phi^{-1}(B)=\tB$, and let $\rG$ be the hyperbolic
geodesic joining the points $\ta$ and $\tb$ in $\Disk$. Then $\fg=\phi(\rG)\subset B$.
\end{lem}

\begin{proof} Suppose, by way of contradiction, that $x\in\rG\sm \tB$.
 Let $C$ be the component of $\ol{\Disk}\sm \ol{\phi^{-1}(\rg)}$ which does not contain $x$.
Choose $a_i\to\ta$ and $b_i\to \tb$ in $S^1\cap C$  and let $H_i$ be the hyperbolic geodesic in $\Disk$
joining the points $a_i$ and $b_i$. Then $\lim H_i=\rG$ and $H_i\cap \tB$ is not connected. This contradiction with
Lemma~\ref{jorg} completes the proof.
\end{proof}

\begin{lem}\label{compactfol}
Suppose that $\{\rG_i\}$ is a sequence of hyperbolic chords in $\Gamma$ and
suppose that $x_i\in \rG_i$ such that $\{x_i\}$ converges to $x\in \mathbb D$.
 Then there is a  unique hyperbolic chord $\rG\in\Gamma$ that
contains $x$. Furthermore, $\lim \rG_i=\ol{\rG}$.
\end{lem}
\begin{proof} We may suppose that the sequence $\{\rG_i\}$ converges to
a hyperbolic chord $\rG$ which contains $x$.
 Let $\rg_i\in\KP$
so that $\phi^{-1}(\rg_i)$ is an open arc which joins the endpoints of $\rG_i$.
By Lemma~\ref{contKP}, $\lim \rg_i=\rg\in\kp$. It follows that $\rG$ is the hyperbolic chord joining the endpoints of $\phi^{-1}(\rg)$. Hence $\rG\in\Gamma$.\end{proof}

So we have used the  family  of
$\kp$ chords in $\tU$ to stratify $\mathbb{D}$ to the family $\Gamma$ of hyperbolic chords.
 By  Lemma~\ref{sameb} for  each $\kp$ chord $\rg\subset \HCH(B\cap \partial\tU)$
 its associated hyperbolic chord $\fg=\phi(\rG)\subset B$. Hence, there is a deformation of $\tU$ that maps $\bigcup \KP$
onto $\bigcup\phi(\Gamma)$, which suggests that  components of $\tU\sm \bigcup \phi(\Gamma)$
 naturally correspond to the interiors of the gaps of
the Kulkarni-Pinkall partition. That this is indeed the case is the substance of the next lemma.

\begin{lem}
There is a $1-1$ correspondence between complementary domains
$Z\subset\mathbb D\setminus \bigcup \Gamma$ and the interiors of  Kulkarni-Pinkall gaps
$\HCH( B\cap K)$. Moreover, for each gap $Z$ of $\Gamma$ there exists a unique
$B\in\B^\infty$ such that $Z$ corresponds to the interior of the $\KP$ gap $\HCH(B\cap K)\cap \tU$
in that $\partial Z\cap \Disk=\bigcup\{G\in\Gamma\mid \rg\in\KP \text{ and }
\rg\subset \partial \HCH(B\cap K)\}$ and $\phi(Z)\subset  B$.
\end{lem}

\begin{proof}Let $\rg$ and $\rh$ be two distinct $\kp$ chords in the boundary of the
gap $\HCH(B\cap K)$ for some $B\in\B^\infty$. Let $\{a,b\}$ and $\{c,d\}$ be the endpoints
of $\phi^{-1}(\rg)$ and $\phi^{-1}(\rh)$, respectively. Then $\rG$ has endpoints $\{a,b\}$
and $\rH$ has endpoints $\{c,d\}$.
There exist disjoint irreducible arcs
$A$ and $C$ in $\partial\Disk$ between the sets $\{a,b\}$ and $\{c,d\}$.
Since $\rg$ and $\rh$ are contained in the same gap, no hyperbolic leaf
of $\Gamma$ has one endpoint in $A$ and the other endpoint in $C$. Hence there exists a gap $Z$
of $\Gamma$ whose boundary includes the hyperbolic chords $\rG$ and $\rH$. It now follows easily that for any
$\rg'\in\KP$ which is contained in the boundary of the same gap $\HCH(B\cap K)$,
$\rG'$ is contained in the boundary of $Z$. Hence the $\KP$ gap $\HCH(B\cap K)$ corresponds to the
gap $Z$ of $\Gamma$. Conversely, if $Z$ is a gap of $\Gamma$ in $\Disk$ then a similar argument,
together with Lemmas~\ref{contKP} and \ref{lenseunion}, implies that $Z$ corresponds to
a unique gap $\HCH(B\cap K)$ for some $B\in\B^\infty$. The rest of the Lemma now follows from Lemma~\ref{sameb}.
\end{proof}

So if $\tU=\sphere\sm K$ is endowed with the hyperbolic metric induced by $\phi$,
then there exists a family of geodesic chords that share the same endpoints as elements of $\KP$.
 The complementary domains
of $\tU\sm \bigcup\{\fg\mid \rg\in\KP\}$ corresponds to the Kulkarni-Pinkall gaps. We summarize the results:

\begin{thm}\label{Hypmain} Suppose that $K\subset\complex$
is a non-separating continuum and let $\tU$ be its
complementary domain in the Riemann sphere. There exists a family $\phi(\Gamma)$ of
hyperbolic chords  in the hyperbolic metric on $\tU$
such that for each $\fg\in\phi(\Gamma)$ there exists $B\in\B^\infty$ and $\rg\subset \HCH(B\cap \partial\tU)$
so that $\fg$ and $\rg$ have the same endpoints and  $\fg\subset B$.  Each domain $Z$ of
$\tU\setminus \phi(\Gamma)$  naturally corresponds to a Kulkarni-Pinkall gap $\HCH(B\cap\partial \tU)$
The bounding hyperbolic chords of $Z$ in $\tU$ correspond to the $\kp$
chords (i.e., chords in $\KP$) of $\HCH(B\cap\partial \tU)$.
\end{thm}

In order to obtain Bell's Euclidean foliation \cite{bell76} we could have modified
the $\KP$ family as follows. Suppose that $B\in\B$.
Instead of replacing a $\kp$ chord  $\rg\in\HCH(B\cap K)$  by a geodesic in the
hyperbolic metric on $\tU$, we could have replaced it by
 a straight line segment; i.e, the geodesic in the Euclidean metric.
Then we would have obtained a family of open straight line segments. In so doing we
would have  replaced
the gaps $\HCH(B\cap\partial \tU)$ by $\ECH(B\cap\partial
\tU)$, which is the way in which Bell originally foliated
$\ECH(K)\setminus K$. We hope that the above argument provides a more transparent proof
of Bell's result. Note that both in the hyperbolic and Euclidean case the elements of the
foliation are not necessarily disjoint (hence we use the word \lq\lq foliate''
rather then \lq\lq partition''). However, in both cases every point of $\tU$
is  contained in either a unique chord or in the interior of a unique gap.

\section{$\kp$ chords and prime ends}We will follow the notation from Section~\ref{secKP}
in the case that $K=T(X)$  where $X$ is a plane continuum.
Here we assume, as in the introduction to this paper,  that
$f:\complex\to\complex$ takes continuum $X$ into $T(X)$ with no
fixed points in $T(X)$, and $X$ is minimal with respect to these
properties. We apply the Kulkarni-Pinkall partition to $\tU=\sphere\sm T(X)$.
Recall that  $\KPP=\{\HCH(B\cap K)\cap \tU\mid B\in\B^\infty\}$ is the Kulkarni Pinkall partition of $\tU$
as given by Theorem~\ref{KPthrm}.

Let $B^\infty\in\B^\infty$ be the  maximal ball such that $\infty\in\HCH(B^\infty\cap K)$.
 As before we use balls on the sphere.
In particular, straight lines in the plane
 correspond to circles on the sphere containing the point at infinity.
The subfamily of $\KPP$ whose elements are of diameter $\leq
\delta$ in the spherical metric is denoted by $\KPP_\delta$. \index{KPPD@$\KPP_\delta$}
 The subfamily
of  chords in $\KP$  of diameter $\leq \delta$ is denoted by $\KP_\delta$. \index{KP chord@$\KP_\delta$}

By Lemma~\ref{contKP} we  know that the families $\KP$ and $\KPP$ have nice continuity
properties. However, $\KP$ and $\KPP$ are not closed in
the hyperspace of compact subsets of $\complex^\infty$: a sequence
of chords or hulls may converge to a point in the boundary
of $\tU$ (in which case it  must be a null sequence).

\begin{prop}[Closedness]\label{compactness}
Let  $\{\rg_i\}$ be a convergent sequence of distinct elements in  $\KP_\delta$,
 then either $\rg_i$ converges to a chord $\rg$ in
$\KP_\delta$ or $\rg_i$ converges to a point of $X$.
In the first case, for  large $i$ and $\delta$ sufficiently small,
$\var(f,\rg,T(X))=\var(f,\rg_i,T(X))$.
\end{prop}
\begin{proof}
By Lemma~\ref{contKP},  we know that the first conclusion holds if
$\rg=\lim \rg_i$ contains a point of $\tU$. Hence we only need to
consider the case when $\lim \rg_i=\rg\subset
\partial \tU\subset T(X)$. If the diameter of $\rg_i$
converged to zero, then $\rg$ is a point as desired. Assume that this
is not the case and let $B_i$ be the maximal ball that contains
$\rg_i$. Under our assumption, the diameters of $\{B_i\}$ do not decay to
zero. Then $\lim B_i=B\in\B^\infty$ and
it follows $\lim \rg_i$ is a piece of a round circle which crosses $\partial B$
perpendicularly.  Hence  $\lim \rg_i\cap \Int(B)\ne\0$, contradicting
the fact that $\rg\subset \partial\tU\subset T(X)$. Note that for $\delta$ sufficiently small,
$\ol{\rg}\cap f(\ol{\rg})=\0$. Hence, $\var(f,\rg,T(X))$ and $\var(f,\rg_i,T(X))$ are defined for all $i$ sufficiently large. Then last statement in the
Lemma follows from stability of variation (see
Section~\ref{compofvar}).
\end{proof}
\begin{cor}\label{small} For each $\e>0$, there exist $\delta>0$
such that for all $\rg\in\KP$ with $\rg\subset B(T(X),\delta)$, $\dm(\rg)<\e$.
\end{cor}
\begin{proof}
Suppose not, then there exist $\e>0$ and  a sequence $\rg_i$ in $\KP$ such that $\lim
\rg_i\subset X$ and $\dm(\rg_i)\geq\e$ a contradiction to Proposition~\ref{compactness}.
\end{proof}

 The proof of the following well-known proposition is omitted.
 \begin{prop}\label{smallcutoff} For each $\e>0$ there exists $\delta>0$
 such that for each open arc $A$ with distinct endpoints $a,b$ such that $\cl{A}\cap T(X)=\{a,b\}$
 and $\dm(A)<\delta$, $T(T(X)\cup A)\subset B(T(X),\e)$.
 \end{prop}

 \begin{prop} \label{smallgeometric} Let $\e$, $\delta$ be as in Proposition~\ref{smallcutoff}
 above with $\delta<\e/2$ and let $B\in\B^\infty$.
Let $A$ be a crosscut of $T(X)$
 such that $\dm(A)<\delta$. If $x\in T(A\cup T(X))\cap \HCH(B\cap T(X))\sm T(X)$
and $d(x,A)\geq\e$,
  then the radius of $B$ is less than $\e$. Hence, $\dm(\HCH(B\cap T(X)))<2\e$.
 \end{prop}
 \begin{proof}
 Let $z$ be the center of $B$.
 If  $d(z,T(X))<\e$ then $\dm(B)<2\e$ and we are done. Hence, we may assume that $d(z,T(X))\ge \e$.
 We will show that this leads to a contradiction. By Proposition~\ref{smallcutoff} and our choice of $\delta$,
 $z\in \sphere\sm T(A\cup X)$.
  The straight line segment $\ell$ from $x$ to $z$ must cross $T(X)\cup A$ at some point $w$.
 Since the segment $\ell$ is in the interior of the maximal ball $B$, it is
 disjoint from $T(X)$, so $w\in A$.
 Hence $d(x,w)\geq\e$ and, since $x\in B$,  $B(w,\epsilon)\subset B$. This is  a contradiction since
 $A\subset B(w,\delta)$ and $\delta<\epsilon/2$ so $\overline A$ would be contained in the interior
 of $B$ which is impossible since $A$ is a crosscut of $T(X)$.
  \end{proof}

\begin{prop}\label{crossing} Let $C$ be a crosscut of $T(X)$
and let $A$ and $B$ be disjoint closed sets in $T(X)$ such that
$\cl{C}\cap A\not=\0\not= \cl{C}\cap B$.
For each $x\in C$,  let $F_x\in\KPP$ so that $x\in F_x$. If each $\ol{F_x}$ intersects $A\cup B$,
then there exists an $F_\infty\in\KPP$ such that $\ol{F_\infty}$ intersects $A$,  $B$
and $\ol{C}$.
\end{prop}
\begin{proof}
Let $a\in A,b\in B$ be the endpoints of $\cl C$. Let $C_a,
C_b\subset C$ be the set of points $x\in C$ such that $\ol{F_x}$ intersects $A$ or
$B$, respectively. Then $C_a$ and $C_b$ are closed subsets by
Proposition~\ref{compactness}. Note that $d(A,B)>0$. If $C_a=\0$, choose $x_i\in C$ converging to $a\in A\cap\ol{C}$.
Then $\ol{F_{x_i}}\cap B_\infty\ne\0$ and $\lim F_{x_i}=\ol{F_\infty}\subset \HCH(B_\infty\cap K)\in\KPP$,  $\lim B_i=B_\infty$ and where $B_\infty,B_i\in\B^\infty$ such that $F_{x_i}\subset\HCH(B_i\cap T(X))$
by Lemma~\ref{contKP}.
Then $\ol{F_\infty} \cap B\ne\0 $ and $a\in A\cap \ol{C}\cap \ol{F_\infty}$.
Suppose now $C_a\ne\0\ne C_b$. Then $C_a$ and $C_b$ are closed and, since $C$ is connected,
$C_a\cap C_b\ne\0$. Let $y\in C_a\cap C_b$. Then $\ol{F_y}\cap A\ne\0\ne \ol{F_y}\cap B$ and $y\in F_y\cap C$.
\end{proof}

Proposition~\ref{crossing} allows us to replace small crosscuts
which essentially cross a prime end $\mc{E}_t$  with non-trivial principal
continuum  by small nearby $\kp$ chords which also essentially cross
$\mc{E}_t$. For if $C$ is a small crosscut of
$T(X)$ with endpoints $a$ and $b$ which crosses the
external ray $R_t$ essentially, let $A$ and $B$ be the closures of
the sets in $T(X)$ accessible from $a$ and $b$, respectively by small
arcs missing $R_t$. If the $F_\infty$ of proposition~\ref{crossing}
is a gap $\HCH(B\cap T(X))$, then a $\kp$ chord in its boundary
crosses $R_t$ essentially.

Fix a Riemann map $\varphi:\disk^\infty\to \tU=\rsphere\sm T(X)$
with $\varphi(\infty)=\infty$. Recall that an external ray $R_t$ is the image
of the radial line segment with argument $2\pi t i$ under the map $\varphi$.

\begin{prop} \label{var0}  Suppose the external ray $R_t$ lands on $x\in T(X)$, and
$\{\rg_i\}_{i=1}^\infty$ is a sequence of crosscuts of $T(X)$ converging to $x$
such that there exists a null sequence of arcs $A_i\subset \complex\setminus T(X)$
joining $\rg_i$
to $R_t$.
  Then for
sufficiently large $i$, $\var(f,\rg_i,T(X))=0$. \end{prop}

\begin{proof}   Since $f$ is fixed point free
on $T(X)$ and $f(x)\in T(X)$, we may choose  a connected neighborhood
$W$ of $x$ such that $f(\cl{W})\cap(\cl{W}\cup R_t)=\0$.  For
sufficiently large $i$, $A_i\cup \rg_i\subset W$.  Then for each such $i$ there exists
 a junction $J_i$ starting from a point in $\rg_i$, staying in $W$ close to $A_i$
until it reaches $R_t$, then following $R_t$ to $\infty$.  By our
choice of $W$, $\var(f,\rg_i,T(X))=0$.
 \end{proof}

\begin{prop}\label{pomm92} Suppose that for an external ray $R_t$ we have $R_t\cap \linebreak[4]
\Int( \ECH(T(X)))
\ne \0$. Then there exists $x\in R_t$ such that the $(T(X),x)$-end of $R_t$ is contained in
$\ECH(T(X))$. In particular there exists a chord $\rg\in\KP$ such that $R_t$ crosses $\rg$
essentially.
\end{prop}

\begin{proof} External rays in $\tU$ correspond to geodesic half-lines starting at
infinity in the hyperbolic metric on $\sphere\sm T(X)$. Half-planes are conformally equivalent to disks.
Therefore, J{\o}rgensen's lemma applies: the intersection of $R_t$
with a halfplane is connected, so it is a half-line. Since the
Euclidean convex hull of $T(X)$ is the intersection of all half-planes containing
$T(X)$, $R_t\cap\ECH(T(X))$ is connected.
\end{proof}

\begin{lem}\label{chordlimit}  Let $\mc{E}_t$ be a channel (that is, a prime end such that
$\pr(\mc{E}_t)$ is non-degenerate) in $T(X)$.
Then for each $x\in\pr(\mc{E}_t)$, for every
$\delta>0$, there is a chain $\{\rg_i\}_{i=1}^\infty$ of chords
defining $\mc{E}_t$ selected from $\KP_\delta$ with  $\rg_i\to
x\in \bd T(X)$.
\end{lem}
\begin{proof} Let $x\in\pr(\mc{E}_t)$ and let $\{C_i\}$ be a defining chain of crosscuts
for $\pr(\mc{E}_t)$ with $\{x\}=\lim C_i$. By
Proposition~\ref{crossing}, in particular by the remark following
the proof of that proposition, there is a sequence $\{\rg_i\}$ of
$\kp$ chords such that $d(\rg_i, C_i)\to 0$ and $\pr(\mc{E}_t)$ crosses
each $\rg_i$ essentially. By Proposition~\ref{smallgeometric}, the
sequence $\rg_i$ converges to $\{x\}$.
\end{proof}

\begin{lem} \label{cutoff} Suppose an external ray $R_t$ lands on $a\in T(X)$ with
$\{a\}=\pr(\mc{E}_t)\not=\im(\mc{E}_t)$.  Suppose
$\{x_i\}_{i=1}^\infty$ is a collection of points in $\tU$
with $x_i\to x\in\im(\mc{E}_t)\sm \{a\}$ and $\phi^{-1}(x_i)\to
t$. Then  there is a sequence of $\kp$ chords
$\{\rg_i\}_{i=1}^\infty$ such that $\rg_i$ separates $x_i$ from
$\infty$, $\rg_i\to a$ and $\phi^{-1}(\rg_i)\to t$ (for sufficiently large $i$).
\end{lem}

\begin{proof} The existence of the chords $\rg_i$ again follows from the remark
following Proposition~\ref{crossing}. It is easy to see that $\lim
\varphi^{-1}(\rg_i)\to t$.
\end{proof}

\subsection{Auxiliary Continua}\label{auxcont}
We use $\kp$ chords to form Carath\'{e}odory loops around the continuum $T(X)$.

\begin{defn} Fix $\delta>0$.  Define the
following collections of chords:\index{KPd chords@$\KP^{\pm}_\delta$}
$$\KP^+_\delta=\{\rg\in\KP_\delta\mid \var(f,\rg,T(X))\ge 0\}$$
$$\KP^-_\delta=\{\rg\in\KP_\delta\mid \var(f,\rg,T(X))\le 0\}$$ To
each collection of chords above, there corresponds an auxiliary
continuum defined as follows: $$T(X)_\delta=T(X\cup(\cup
\KP_\delta))$$ \index{TXd@$T(X)_\delta$} $$T(X)^+_\delta=T(X\cup(\cup \KP^+_\delta))$$
$$T(X)^-_\delta=T(X\cup(\cup \KP^-_\delta))$$
\end{defn}

\begin{prop} \label{boundary}
Let $Z\in\{T(X)_\delta,T(X)^+_\delta,T(X)^-_\delta\}$, and correspondingly
$\mc{W}\in\{\KP_\delta,\KP^+_\delta,\KP^-_\delta\}$. Then
the following hold: \begin{enumerate} \item $Z$ is a nonseparating
plane continuum. \item $\bd Z\subset T(X)\cup(\cup\mc{W})$. \item Every
accessible point $y$ in $\bd Z$ is either a point of $T(X)$ or a point
interior to a chord $\rg\in\mc{W}$.
\item  If $y\in\partial Z\cap\rg$  with $\rg\in\mc{W}$, then $y$ is accessible,
$\rg\subset \partial Z$
and $\partial Z$ is locally connected at each point of $\rg$.  Hence, if $\varphi:\disk^\infty\to\rsphere\sm Z$
is the Riemann map and $R_t$ is an external ray landing at $y$, then $\varphi$ extends
continuously to an open interval in $S^1$ containing $t$.
\end{enumerate}
\end{prop}

\begin{proof} By Proposition~\ref{compactness}, $T(X)\cup(\cup \mc{W})$ is
compact.  Moreover, $T(X)\cup(\cup \mc{W})$  is connected since each crosscut
$A\in\mc{W}$ has endpoints in $T(X)$.  Hence, the topological hull
$T(T(X)\cup(\cup \mc{W}))$ is a nonseparating plane continuum,
establishing (1).

Since $Z$ is the topological hull of $T(X)\cup(\cup \mc{W})$, no
boundary points can be in complementary domains of $T(X)\cup(\cup
\mc{W})$.  Hence, $\bd Z\subset T(X)\cup(\cup \mc{W})$, establishing
(2). Conclusion (3) follows immediately.

Suppose $y\in\partial Z\cap\rg$ with $\rg\in\mc{W}$. Then $\Sh(\rg)\subset Z$  and there exists
$y_i\in\complex\sm Z$ such that $\lim y_i=y$. We may assume that all the points $y_i$ are
on the ``same side'' of the arc $\rg$  (i.e., $y_i\in\complex\sm \Sh(\rg)$).  This side of $\rg$ is either (1) a limit of $\kp$ chords $\rg_j$, or (2)
there exists a gap $\HCH(B\cap X)$ on this side  with $\rg$ in its boundary. In case (1),
$\rg\subset \Sh(\rg_j)$ and, since $y_i\in \complex\sm Z$ for all $i$, $\rg_j\not\in \mc{W}$.
Hence each $\rg_j\subset \complex\sm Z$ for all $j$. It follows that every point of $\rg$ is accessible, $\rg\subset \partial Z$ and $\partial Z$ is locally connected at each point of $\rg$.
In case (2) there exists a chord $ \rg'\ne\rg$ in the boundary of $\HCH(B\cap X)$ which separates
$\rg$ from infinity. Then $\rg'\not\in\mc{W}$ and  the interior of $\HCH(B\cap X)\subset \complex\sm Z$. Hence the same conclusion follows.

 The last part of (4)
follows from  the proof of Carath\'eodory's theorem
(see \cite{pomm92}).

\end{proof}

\begin{prop}\label{LC} $T(X)_\delta$ is locally connected;
hence, $\bd T(X)_\delta$ is a Carath\'eodory loop.
\end{prop}

\begin{proof}
Suppose that $T(X)_\delta$ is not locally connected. Then
$T(X)_\delta$ has
a non-trivial impression and there exist $0<\e<\delta/2$ and  a   chain $A_i$ of crosscuts of
$T(X)_\delta$ such that $\dm( \Sh(A_i))>\e$ for all $i$. We may assume that $\lim A_i=y\in T(X)_\delta$.

 By Proposition~\ref{boundary} (4)
  we may assume $y\in X$.
Choose $z_i\in\Sh(A_i)$ such that $ d(z_i,y)>\e$.
We can enlarge the crosscut $A_i$ of $T(X)_\delta$ to a crosscut $C_i$ of $T(X)$ as follows.
Suppose that $A_i$ joins the points $a^+_i$ and $a^-_i$ in $T(X)_\delta$. If $a^+_i\in T(X)$, put
$y^+_i=a^+_i$. Otherwise $a^+_i$ is contained in a chord $\rg^+_i\in\KP_\delta$, with endpoints in $T(X)$,
 which is contained in
$T(X)_\delta$. Since $\lim A_i=y$,
  we can select one of these endpoints and call it $y^+_i$ such that
$d(y^+_i,a^+_i)\to 0$. Define $\rg^-_i$ and $y^-_i$ similarly.
Then $\rg^+_i\cup A_i\cup \rg^-_i$ contains a crosscut $C_i$ of $T(X)$ joining the points $y^+_i$ and $y^-_i$
such that $\lim C_i=y$. We claim that $z_i\in\Sh(C_i)$. To see this note that, since $z_i\in Sh(A_i)$,
there exists a halfray $R_i\subset \complex\setminus T(X)_\delta$ joining $z_i$ to infinity
such that $|R_i\cap A_i|$ is an odd number and each intersection is transverse. Since $R_i\cap C_i=R_i\cap A_i$
it follows that $z_i\in\Sh(C_i)$.  Let $\HCH(B_i\cap X)$ be the unique hull of the Kulkarni-Pinkall partition
$\KPP$ which contains $z_i$. Since $\dm(C_i)\to 0$ and $d(z_i,y)>\e$, it follows from
Proposition~\ref{smallgeometric} that $\dm(\HCH(B_i\cap X))<2\e<\delta$. This contradicts the fact that
$z_i\in\complex\setminus T(X)_\delta$ and completes the proof.

\end{proof}

\section{Outchannels} \label{out} Suppose that $f:\complex\to \complex$ is a perfect map,
$X$ is a continuum, $f$ has no fixed point in $T(X)$ and
$X$ which is minimal with respect to $f(X)\subset T(X)$. Fix $\eta>0$ such
that for each $\kp$ chord $\rg\subset T(X)_\eta$, $\ol{g}\cap f(\ol{\rg})=\0$ and $f$ is
 fixed point free on $T(X)_\eta$. In this case we will say that $\eta$ \emph{defines variation near}
\index{defines variation near $X$} $X$ and that the triple
$(f,X,\eta)$ satisfy the \emph{standing hypothesis}. \index{standing hypothesis}\index{fXeta@$(f,X,\eta)$}
In this section we will show that $X$ has at least one
\textit{negative outchannel}, which is defined as follows. Note that for
 each $\kp$ chord $\rg$ in $T(X)_\eta$, $\var(f,\rg,T(X))=\var(f,\rg)$ is defined.

\begin{defn} [Outchannel]  Suppose that $(X,f,\eta)$ satisfy the standing hypothesis.
An {\em outchannel} \index{outchannel}
of the nonseparating plane continuum
$T(X)$ is a prime end $\mc{E}_t$ of $\tU=\rsphere\sm T(X)$ such
that for some chain $\{\rg_i\}$ of crosscuts defining $\mc{E}_t$,
$\var(f,\rg_i,T(X))\not= 0$ for every $i$. We call an outchannel
$\mc{E}_t$ of $T(X)$ a {\em geometric outchannel} \index{outchannel!geometric} iff for
sufficiently small $\delta$, every chord in $\KP_\delta$, which
crosses  $\mc{E}_t$ essentially, has nonzero variation. We call a
geometric outchannel {\em negative} (respectively, {\em positive})
(starting at $\rg\in\kp$)
iff every $\kp$ chord  $\rh \subset T(X)_\eta \cap \ol{\Sh(\rg)}$ , which crosses $\mc{E}_t$
essentially, has negative (respectively, positive) variation.
\end{defn}

\begin{lem} \label{localarcs} Suppose that $(f,X,\eta)$ satisfy the standing hypothesis and $\delta\le\eta$.
Let
$Z\in\{T(X)^+_\delta,T(X)^-_\delta\}$. Fix a Riemann map
$\varphi:\disk^\infty\to\rsphere\sm Z$ such that $\varphi(\infty)=\infty$.
Suppose $R_t$ lands at $x\in\bd Z$. Then there is an open interval
$M\subset\bd\disk^\infty$ containing $t$ such that $\varphi$ can be
extended  continuously over $M$. \end{lem}

\begin{proof}

Suppose that $Z=T(X)_\delta^-$ and $R_t$ lands on $x\in\partial Z$.  By proposition~\ref{boundary}
we may assume that $x\in X$.  Note first that the family of chords in $\kp^-_\delta$ form
a closed subset of the hyperspace of  $\complex\sm X$, by Proposition~\ref{compactness}.
By symmetry, it suffices to show that we can extend $\psi$ over an interval $[t',t]\subset S^1$ for $t'<t$.

Let $\phi:\disk^\infty\to \complex\sm T(X)$ be the Riemann map for $T(X)$. Then there
exists $s\in S^1$ so that the external ray $R_s$ of $\complex\sm T(X)$  lands at $x$. Suppose
first that there exists a chord $\rg\in \kp^-_\delta$ such that $G=\vp^{-1}(\rg)$ has endpoints
 $s'$ and $s$ with $s'<s$. Since $\kp^-_\delta$ is closed,
there exists a minimal $s"\le s'<s$ such that there exists a chord $\rh\in\kp^-_\delta$ so that
$H=\vp^{-1}(\rh)$ has endpoints $s"$ and $s$. Then $\rh\subset \partial Z$ and $\phi$ can be
 extended over an interval $[t',t]$ for some $t'<t$, by Proposition~\ref{boundary} (4).

Suppose next that no such chord $\rg$ exists. Choose a junction $J_x$ for $T(X)^-_\delta$ and
a neighborhood $W$ of $x$ such that $f(W)\cap [W\cup J_x]=\0$.
We will first show that there exists $\nu\le \delta$ such that
$x\in \partial T(X)_\nu$. We may assume that $\nu$ is so small that any chord of $\kp_\nu$ with
endpoint $x$ is contained in $W$.  For suppose that this is not the case. Then
there exists a sequence $\rg_i\in\KP$ of chords such that $x\in\Sh(\rg_{i+1})\subset\Sh(\rg_i)$,
$\lim \rg_i=x$ and $\var(f,\rg_i)>0$ for all $i$. This contradicts Proposition~\ref{var0}.
Hence $x\in\partial T(X)_\nu$ for some $\nu>0$.

By Proposition~\ref{LC}, the boundary of $T(X)_\nu$ is a simple closed curve $S$ which must contain $x$. If there exists a chord $\rh\in\kp_\nu$  with endpoint $x$ such that $H$  has endpoints
$s'$ and $s$ with $s'<s$ then,  since $\rh\subset W$, $f(\rh)\cap J_x=\0$, $\var(f,\rh)=0$ and $\rh\in \kp^-_\delta$, a contradiction. Hence, chords $\rh$  close to $x$ in $S$ so that $H$ has endpoints
less than $s$ are contained in $W$ and, $\var(f,\rh)=0$ by Proposition~\ref{var0}. Hence
a small interval $[x',x]\subset S$, in the counterclockwise order on $S$ is contained in $T(X)^-_\nu$.
It now follows easily that a similar arc exists in the boundary of $T(X)^-_\delta$ and the desired result follows.
\end{proof}

By a \emph{narrow strip} we mean the image of an embedding $h:\{(x,y)\in\Complex\mid x\ge 0 \text{ and } -1\le y\le 1\}\to \Complex$
such that $\lim_{x\to\infty}\dm (h(\{x\}\times [-1,1]))=0$.

\begin{figure}

\includegraphics{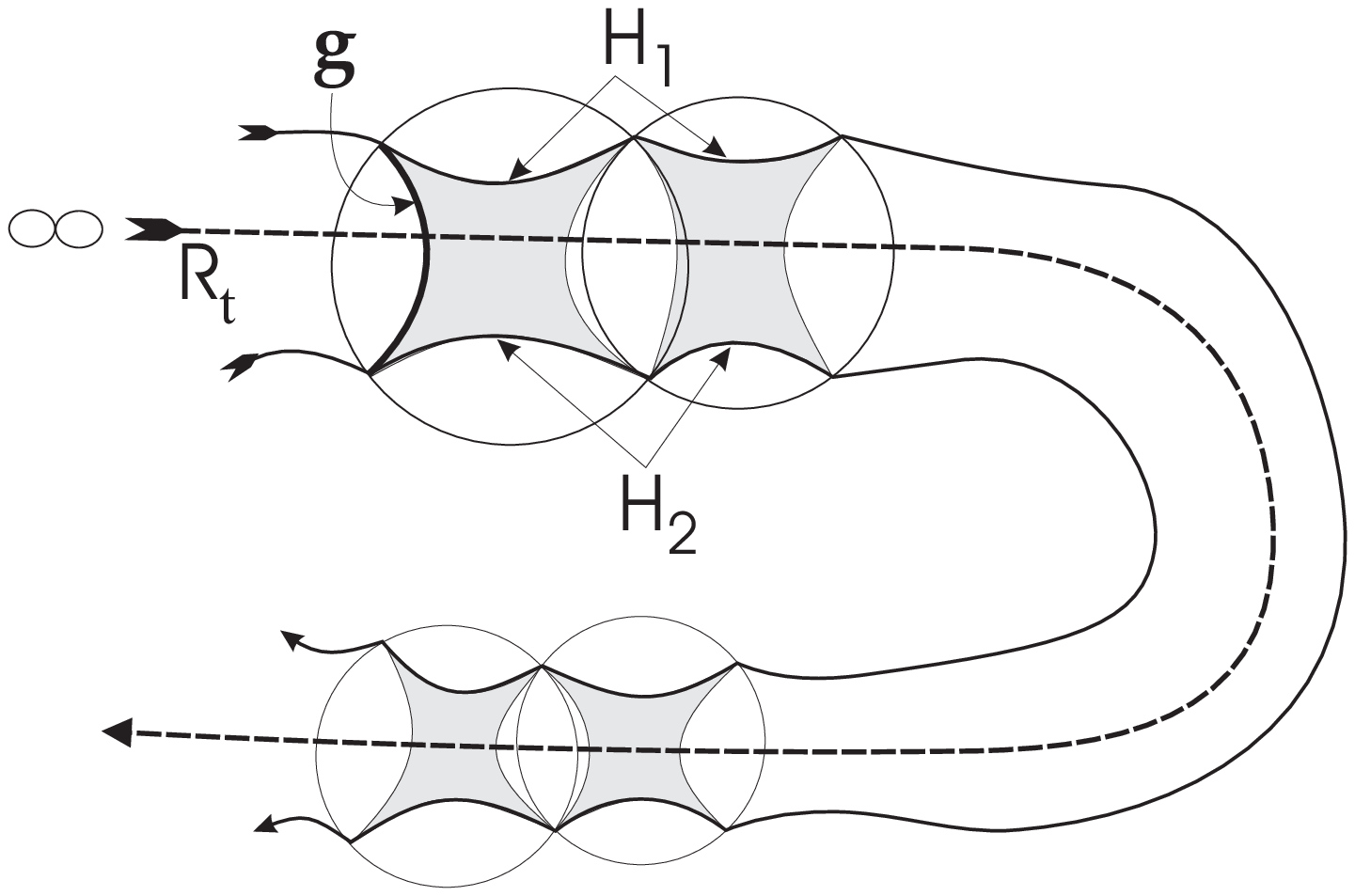}

\caption{The strip $\fS$ from Lemma~\ref{geomoutchannel}} \label{strip}

\end{figure}
\begin{lem} \label{geomoutchannel}Suppose that $(f,X,\eta)$ satisfy the standing hypothesis.  If there is a chord $\rg$ of $T(X)$ of
negative (respectively, positive) variation, such that there is no
fixed point in $T(T(X)\cup \rg)$, then there is a negative
(respectively, positive) geometric outchannel $\mc{E}_t$ of $T(X)$
starting at $\rg$.

Moreover, if $\mc{E}_t$ is a positive (negative) geometric outchannel starting at the $\kp$ chord
$\rg$, and $\fS=\bigcup \{\HCH(B\cap T(X)) \mid \HCH(B\cap T(X))\subset T(X)_\eta\cap \ol{\Sh(\rg)} \text{ and a chord in } \HCH(B\cap T(X)) \text{ crosses } R_t \text{ essentially}\}.$ Then $\fS$ is an infinite narrow strip in the plane whose remainder is contained in $T(X)$ and which
is bordered by a $\kp$ chord and  two halflines $H_1$ and $H_2$ (see figure~\ref{strip}).
\end{lem}

\begin{proof}  Without loss of generality, assume $\var(f,\rg,T(X))=\var(f,\rg)<0$. If $\rg\subset
T(X)_\eta$, put $\rg'=\rg$. Otherwise consider the boundary of
$T(X)_\eta$ which is locally connected by Proposition~\ref{LC} and, hence, a Carath\'eodory loop. Then a continuous extension $g:S^1\to\partial T(X)_\eta$    of the Riemann map  $\phi:\disk^\infty\to \rsphere\sm T(X)_\eta$ exists. Whence the boundary of
$T(X)_\eta$ contains a sub-path $A=g([a,b])$,  which is contained in $\ol{\Sh(\rg)}$, whose endpoints coincide with the endpoints of $\rg$. Note that for each component $C$ of $A\sm X$,
$\var(f,C)$ is defined. Then it follows from Proposition~\ref{crossunder}, applied to a
Carath\'eodory path, that there exists a component $C=\rg'$ such that $\var(f, \rg')<0$. Note that
$\rg'$ is a $\kp$ chord contained in the boundary of $T(X)_\eta$.

To see that a geometric outchannel, starting with $\rg$ exists,
note that for any chord $\rg''\subset T(X)_\delta$ with $\var(f,\rg'',X)<0$, if $\rg''=\lim \rg_i$,
then there exists $i$ such that for any chord $\rh$ which separates $\rg_i$ and $\rg''$
in $U^\infty$, $\var(f,
\rh,X)<0$ (see Proposition~\ref{var0}). Also by the argument above, if
in addition  $\rg''\subset \HCH(B\cap T(X))$, where $\HCH(B\cap T(X))$ is a gap, such that $\rg''$ separates
$\HCH(B\cap T(X))\sm\rg''$ from infinity in $U^\infty$, then there exists $\rh\ne\rg''$
in $\HCH(B\cap T(X))$ such that $\rg''$ separates $\rh$ from infinity in $U^\infty$ and
$\var(f,\rh,X)<0$.  The remaining conclusions of the Lemma  follow from these two facts.
\end{proof}

\subsection{Invariant Channel in $X$}
We are now in a position to prove Bell's principal result on any
possible counter-example to the fixed point property, under our
standing hypothesis.

\begin{lem} \label{invchannel}
Suppose $\mc{E}_t$ is a geometric outchannel of $T(X)$ under
$f$. Then the principal continuum $\pr(\mc{E}_t)$ of $\mc{E}_t$ is
invariant under $f$. So $\pr(\mc{E}_t)=X$. \end{lem}

\begin{proof}
Let $x\in\pr(\mc{E}_t)$. Then for some chain $\{\rg_i\}_{i=1}^\infty$
of crosscuts defining $\mc{E}_t$ selected from $\KP_\delta$, we
may suppose $\rg_i\to x\in \bd T(X)$ (by Lemma~\ref{chordlimit}) and $\var(f,\rg_i,X)\not=0$ for
each $i$. The external ray $R_t$ meets all $\rg_i$ and there is, for each $i$, a junction
 from $\rg_i$  which ``parallels" $R_t$. Since
$\var(f,\rg_i,X)\not=0$, each $f(\rg_i)$ intersects $R_t$. Since
$\dm(f(\rg_i))\to 0$, we have $f(\rg_i)\to f(x)$ and
$f(x)\in\pr(\mc{E}_t)$. We conclude that $\pr({\mc{E}_t})$ is
invariant.\end{proof}

\begin{thm} [Dense channel, Bell] \label{densechannel}\index{channel!dense}
Suppose that $(X,f,\eta)$ satisfy  our standing hypothesis. Then $T(X)$ contains a negative
geometric outchannel; hence, $\bd \tU=\bd T(X)=X=f(X)$ is an
indecomposable continuum.
\end{thm}

\begin{proof}
  By Lemma~\ref{LC} $\bd
T(X)_\eta$ is a Carath\'eodory loop. Since $f$ is fixed point free on
$T(X)_\eta$, $\ind(f,\bd T(X)_\eta)=0$. Consequently, by
Theorem~\ref{I=V+1} for Carath\'eodory loops, $\var(f,\bd
T(X)_\eta)=-1$. By the summability of variation on $\bd T(X)_\eta$, it
follows that on some chord $\rg\subset \bd T(X)_\eta$, $\var(f,\rg,T(X))
<0$. By Lemma~\ref{geomoutchannel}, there is a negative geometric
outchannel $\mc{E}_t$ starting at $\rg$.

Since $\pr(\mc{E}_t)$ is invariant under $f$ by
Lemma~\ref{invchannel}, it follows that $\pr(\mc{E}_t)$ is an
invariant subcontinuum of $\bd \tU\subset\bd T(X)\subset X$. So by
the minimality condition in our Standing Hypothesis, $\pr(\mc{E}_t)$
is dense in $\bd \tU$. Hence, $\bd \tU=\bd T(X)=X$ and
$\pr(\mc{E}_t)$ is dense in $X$.  It then follows from a theorem of
Rutt \cite{rutt35} that $X$ is an indecomposable continuum.
\end{proof}

\begin{thm}\label{nicebd} Suppose that $(X,f,\eta)$ satisfy  our standing Hypothesis and $\delta\le \eta$. Then the boundary of $T(X)_\delta$ is a simple closed curve.
The set of accessible points in the boundary of each of
$T(X)^+_\delta$ and $T(X)^-_\delta$ is an at most countable union of continuous
one-to-one images of $\real$.
\end{thm}

\begin{proof} By Theorem~\ref{densechannel}, $X$ is indecomposable, so it has no
cut points. By Proposition~\ref{LC}, $\bd T(X)_\delta$ is a
Carath\'eodory loop. Since $X$ has no cut points, neither does
$T(X)_\delta$.  A Carath\'eodory loop without cut points is a simple
closed curve.

Let $Z\in\{T(X)_\delta^+,T(X)_\delta^-\}$ with $\delta\le \eta$.  Fix a Riemann map
$\phi:\disk^\infty\to\rsphere\sm Z$ such that $\phi(\infty)=\infty$.
Corresponding to the choice of $Z$, let
$\mc{W}\in\{\KP^+_\delta,\KP^-_\delta\}$.  Apply
Lemma~\ref{localarcs} and find the maximal collection $\mc{J}$ of
disjoint open subarcs of $\bd\disk^\infty$ over which $\phi$ can be
extended continuously. The collection $\mc{J}$ is countable. Since
$X$ has no cutpoints the extension is one-to-one over $\cup\mc J$.
Since angles that correspond to accessible points are dense  in
$\bd\disk^\infty$, so is $\cup\mc{J}$. If $Z=T(X)_\delta^+$, then it is
possible that $\cup\mc{J}$ is all of $\bd\disk^\infty$ except one
point, but it cannot be all of $\bd\disk^\infty$ since there is at
least one negative geometric outchannel by
Theorem~\ref{densechannel}.
\end{proof}

Theorem~\ref{nicebd} still leaves open the  possibility that
$Z\in\{T(X)_\delta^+,T(X)_\delta^-\}$ has a very complicated boundary. The
set $C=\bd\disk^\infty\setminus\cup\mc J$ is compact and
zero-dimensional. Note that $\phi$ is discontinuous at points in
$C$.   We may call $C$ the set of outchannels of $Z$. In principle,
there could be an uncountable set of outchannels, each dense in $X$.
The one-to-one continuous images of half lines in  $\real$ lying in $\bd Z$ are the
``sides" of the outchannels. If two elements $J_1$ and $J_2$ of the
collection $\mc{J}$ happen to share a common endpoint $t$, then the
prime end $\mc{E}_t$ is an outchannel in $Z$, dense in $X$, with
images of half lines
$\phi(J_1)$ and $\phi(J_2)$ as its sides.  It seems possible that an
endpoint $t$ of $J\in\mc{J}$ might have a sequence of elements $J_i$
from $\mc{J}$ converging to it.  Then the outchannel $\mc{E}_t$
would have only one (continuous) ``side." Such exotic possibilities
are eliminated in the next section.

In the lemma below we summarize several of the results in this section and show that an arc component $K$  of the set of accessible points of the boundary of
$T(X)_{\delta}^-$ is efficient in connecting close points in $K$.

\begin{prop}\label{smallarc}Suppose that $(X,f,\eta)$ satisfy  our standing Hypothesis,
 that the boundary of $T(X)_\delta^-$ is not a simple closed curve, $\delta\le\eta$ and that
$K$ is an  arc component of the boundary of $T(X)_\delta^-$ so that $K$ contains an accessible point. Let $\varphi:\disk^\infty\to \rsphere\sm T(X)_\delta^-$ be a conformal map such that $\varphi(\infty)=\infty$. Then:
\begin{enumerate}
\item \label{extpsi}  $\varphi$ extends continuously and injectively to a map
$\tilde{\varphi}:\tilde{\disk}^{\infty}\to \tilde{U}^{\infty}$,
where $\tilde{\disk}^{\infty}\sm \disk^\infty$ is a dense and open
subset of $S^1$ which contains $K$ in its image.  Let $\tilde{\varphi}^{-1}(K)=(t',t)\subset S^1$
with $t'<t$ in the counterclockwise order on $S^1$. Hence $\tilde{\varphi}$ induces an order
$<$ on $K$. If $x<y\in K$, we denote by $<x,y>$ the subarc of $K$ from $x$ to $y$
and by $<x,\infty>=\cup_{y>x} <x,y>$.
\item $\mc{E}_t$ and $\mc{E}_{t'}$ are positive geometric outchannels of $T(X)$.
\item Let $R_t$ be the external ray of $T(X)_\delta^-$ with argument $t$.  There exists $s\in R_t$,
$B\in\B^\infty$ and $\rg\in\kp$ such that $s\in \rg\subset \HCH(B\cap X)$ and
$s$ is the last point of $R_t$ in $\HCH(B\cap X)$ (from $\infty$), $\rg$ crosses $R_t$ essentially and for each $B'
\in \B$ with $\HCH(B'\cap X)\sm X\subset \Sh(\rg)$, $\dm(B')<\delta$.
\item\label{conhul} There exists $\hx\in K$ such that if $B'\in\B^\infty$ with $\Int(B')\subset \Sh(\rg)$, then
$\HCH(B'\cap X)\cap <\hx,\infty>$ is a compact ordered subset of $K$ so that if $C$ is a component of $<\hx,\infty> \sm \, \HCH(B'\cap X)$ with two endpoints in $\HCH(B'\cap X)$,  $C\in \KP_\delta^-$.
\item \label{esscross}
 Let $\B^\infty_t\subset\B^\infty$ be the collection of all $B\in\B^\infty$ such that
 $ R_t$ crosses a chord in the boundary of $\HCH(B\cap X)$ essentially and $\Int(B)\subset\Sh(\rg)$.
 Then $\fS=\bigcup_{B\in \B^\infty_t}  \HCH(B\cap X)$ is a narrow strip
 \index{narrow strip} in the plane, bordered by two halflines $H_1$ and $H_2$,
  which compactifies on $X$ and one of $H_1$ or $H_2$ contains the set $<\hx',\infty>$ for some $\hx'\in K$.\\
 In particular,  if $\max(\hx,\hat{x}')<p<q$ and $\dm(<p,q>) >2\delta$,
 then there exists a chord $\rg\in\KP$ such that one endpoint of $\rg$ is in $<p,q>$ and
 $\rg$  crosses $R_t$ essentially.
\end{enumerate}
The same conclusion holds for $T(X)_\delta^+$ since its boundary cannot be a simple closed curve.
\end{prop}

\begin{proof}By Proposition~\ref{boundary} and Theorem~\ref{nicebd}, and its
proof, $\varphi$ extends continuously and injectively to a map
$\tilde{\varphi}:\tilde{\disk}^{\infty}\to \tilde{U}^{\infty}$ and (\ref{extpsi}) holds.

By Lemma~\ref{localarcs}, the external ray $R_t$ does not land. Hence there exist
a  chain $\rg_i$ of $\KP_\delta$ chords which define the prime end $\mc{E}_t$.
If for any $i$
 $\var(f,\rg_i)\le 0$,
then $\rg_i\subset T(X)_\delta^-$ a contradiction with the definition of $t$.
Hence $\var(f,\rg_i)>0$ for all $i$ sufficiently small and $\mc{E}_t$ is a positive geometric
outchannel by the proof of Lemma~\ref{geomoutchannel}. Hence (2) holds.

Since $R_t$ does not land, we can choose $s'\in R_t$ and $B\in \B^\infty$
such that $s'\in\HCH(B\cap X)$, a chord $\rg'\subset\HCH(B\cap X)$ crosses $R_t$ essentially and  if $B'\in\B^\infty$ such that $\HCH(B'\cap X)\subset \Sh(\rg')$  then $\dm(B')<\delta/2$.   Let $s\in\rg\subset\HCH(B\cap X)$ be the last point on $R_t\cap \HCH(B\cap X)$, starting from $\infty$.  Then (3) holds.

Let  $B'\in\B^\infty$ with $\Int(B')\subset\Sh(\rg)$, $\hx\not\in B'$ and $<\hx,\infty>\cap \HCH(B'\cap X)$ not connected.  Let $\rh\in\kp$ with $\rh\subset \HCH(B'\cap X)$ which separates
$\HCH(B'\cap X)\sm \rh$ from $\infty$ in $U^\infty$. Let $a$ and $b$ be the endpoints of $\rh$.
Then $\{a,b\}\subset <\hx,\infty>$. If there is a chord $\rh'\in\kp$ with endpoints $a$ and $b$
and $\rh'\subset <\hx,\infty>$ we are done. If such a $\kp$ chord does not exist, then
$\var(f,\rh)>0$. By Lemma~\ref{geomoutchannel}, there is a geometric outchannel $\mc{E}_{t"}$
starting at $\rh$. This outchannel disconnects the arc $<a,b>$ between $a$ and $b$,
a contradiction. Hence (4) holds.

 Next choose $\hat{x}'\in K$ such that each point of  $<\hat{x}',\infty>$ is accessible from $\Sh(\rg)$. Then each subarc $<p,q>$ of $<\hat{x},\infty>$ of diameter bigger than $2\delta$ cannot be contained in a single element of the $\KPP$ partition.
Hence there exists a $\kp$ chord $\rg$ which crosses $R_t$ essentially and has one endpoint in
$<p,q>$.

Note that for each chord $\rh\subset \Sh(\rg)$ which crosses $R_t$ essentially, $\var(f,\rh)>0$.
By Lemma~\ref{geomoutchannel},
$\bigcup_{B\in \B^\infty_t}  \HCH(B\cap X)$ is a  strip in the plane, bordered by two halflines $H_1$, $H_2$, which compactify on $X$. These two halflines, consist of chords in $\kp_\delta$ and points in $X$,   one of which, say $H_1$ meets
$<\hx',\infty>$.  If $<\hx',\infty>$ is not contained in $H_1$ then, as in the proof of (4),
there exists a chord $\rh\subset P_1$ with $\var(f,\rh)>0$ joining two points of $x,y\in<\hx,\infty>$ with
all the points of $<x,y>$ accessible from $\Sh(\rh)$.  As above this leads to the contradiction that
$<x,y>$ contains a chord of positive variation and the proof is complete.
\end{proof}

\section{Uniqueness of the Outchannel}
Theorem~\ref{densechannel} asserts the existence of at least one
negative geometric outchannel which is dense in $X$.  We show
below that there is exactly one geometric outchannel, and that its
variation is $-1$.  Of course, $X$ could have other dense
channels, but they are ``neutral" as far as variation is
concerned.

\begin{thm} [Unique Outchannel] \label{outchannel}\index{outchannel!uniqueness}
Suppose that $(X,f,\eta)$ satisfy the standing hypothesis. Then there exists a
unique geometric outchannel $\mc{E}_t$ for $X$, which is dense in
$X=\bd T(X)$. Moreover, for any sufficiently small chord $\rg$ in any
chain defining $\mc{E}_t$, $\var(f,\rg,X)=-1$, and for any
sufficiently small chord $\rg'$ not crossing $R_t$ essentially,
$\var(f,\rg',X)=0$.
\end{thm}

\begin{figure}

\includegraphics{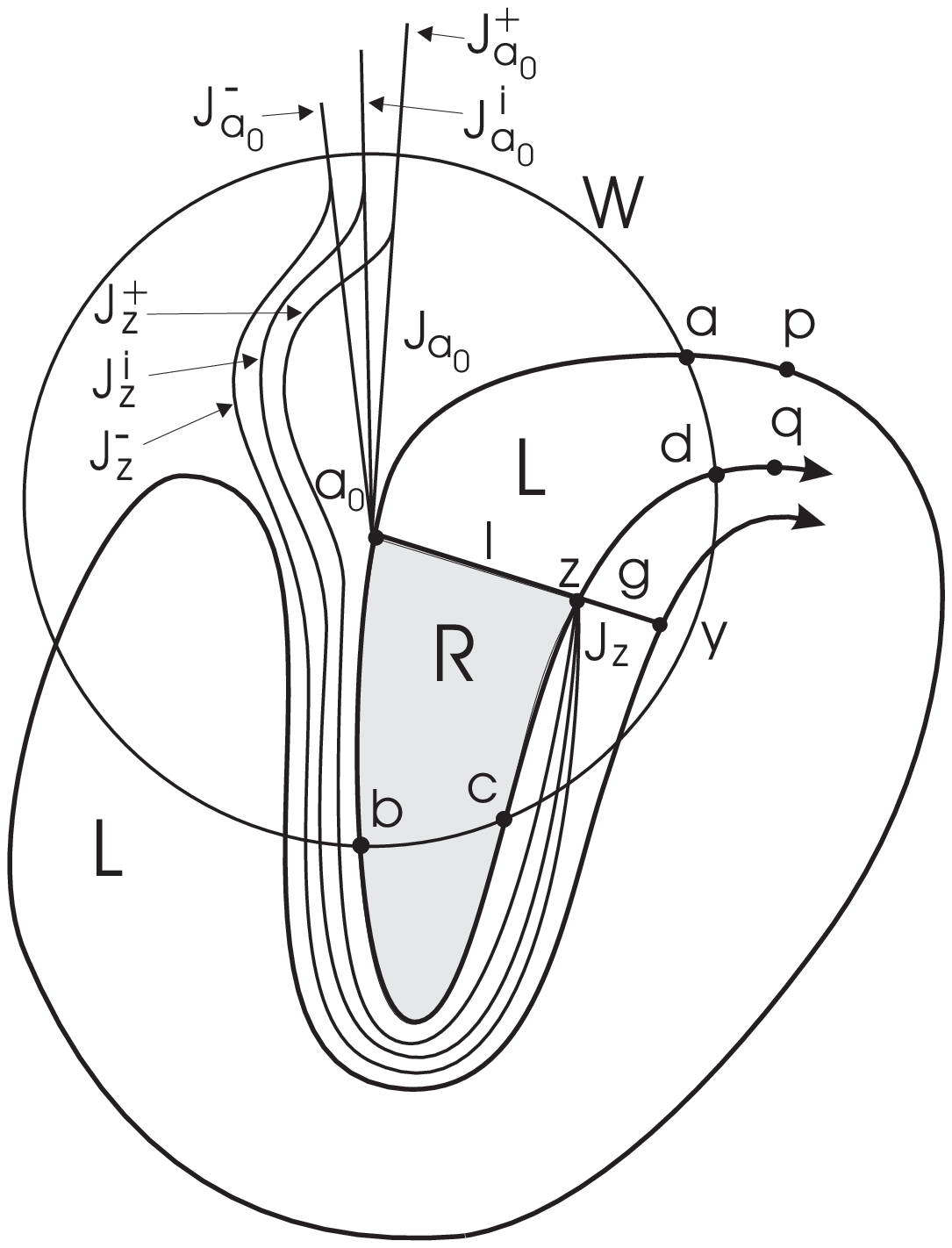}

\caption{Uniqueness of the negative outchannel.} \label{outpic}

\end{figure}

\begin{proof}
Suppose by way of contradiction that $X$ has a positive outchannel.
Let $0<\delta\le\eta$ such that  if $M\subset T(B(T(X),2\delta))$ with
$\dm(M)<2\delta$, then $f(M)\cap M=\0$. Since $X$ has a positive
outchannel, $\bd T(X)_{\delta}^-$ is not a simple closed curve. By
Theorem~\ref{nicebd} $\bd T(X)_{\delta}^-$ contains an arc component
$K$ which is the one-to-one continuous image of $\real$. Note that
each point of $K$ is accessible.

Let $\varphi:\disk^\infty\to V^\infty=\complex\setminus
T(X)_{\delta}^-$ a conformal map.  By Proposition~\ref{smallarc}, $\varphi$ extends continuously and injectively to a map
$\tilde{\varphi}:\tilde{\disk}^{\infty}\to \tilde{V}^{\infty}$,
where $\tilde{\disk}^{\infty}\sm \disk^\infty$ is a dense and open
subset of $S^1$ which contains $K$ in its image. Then
$\tilde{\varphi}^{-1}(K)=(t',t)\subset S^1$ is an open arc with
$t'<t$ in the counterclockwise order on $S^1$ (it could be that $(t',t)=S^1\sm\{t\}$ and $t=t'$). By abuse of notation, let $<$ denote the
order in $K$ induced by $\tilde{\varphi}$ and for $x<y$ in $K$,
denote the arc in $K$ with endpoints $x$ and $y$ by $\langle
x,y\rangle$. For $x\in K$, let $\langle x,\infty\rangle =\cup_{y>x} \langle
x,y\rangle $

Let $\mc{E}_t$ be the prime-end corresponding to $t$. By Proposition~\ref{smallarc},
$\pr(\mc{E}_t)$ is a positive geometric outchannel and, hence, by
Lemma~\ref{invchannel}, $\pr(\mc{E}_t)=X$. Let
$R_t=\varphi(re^{it})$, $r>1$,  be the external conformal ray
corresponding to the prime-end $\mc{E}_t$ of $T(X)^-_\delta$. Since $\cl{R_t}\sm R_t=X$
and  the small chords $\rg_x$ which define $\pr({\mc{E}}_t)$ have at least one endpoint in $K$, cross $R_t$
essentially at $x$ and have diameter going to zero as $x$ approaches $X$ along $R_t$
(by Proposition~\ref{smallgeometric} and Lemma~\ref{chordlimit}).
Hence, $X=\ol{R_t}\sm R_t=\linebreak \ol{<x,\infty>}\ \sm <x,\infty>$.

By Proposition~\ref{smallarc}, there is  $s\in R_t$ such that if $B\in\B^\infty$ such that
 $\HCH(B\cap X)\cap [(X,s)$-end of $R_t]\ne\0$,
then $\dm(B)<\delta/2$.  Let $\B^\infty_t=\{B\in\B^\infty\mid \HCH(B\cap X) \text{ contains
a chord } \rg \text{ such that } \rg \text{ crosses the }\linebreak (X,s) \text{-end of }
 R_t
 \text{ essentially}\}$.

 By Proposition~\ref{smallarc} there exists  $\hx\in K\cap X$ such that
 for each arc $A\subset <\hx,\infty>$ with diameter $>2\delta$, there is a $\kp$ chord
 $\rg$ which contains a point of $A$ as an endpoint and crosses $R_t$ essentially.

Let $a_0\in K\cap X$ so that $a_0>\hx$ and $J_{a_0}$ is a junction of $T(X)_\delta^-$.
Let $W$ be an open disk, with simple closed curve boundary, about $a_0$ such that
$\dm(W)<\delta/4$ and $f(\ol{W})\cap [\ol{W}\cup J_{a_0}]=\0$. Let $a<a_0<b$ in $K\cap\partial W$ such that $<a,b>$ is the component of $K\cap\ol{W}$ which contains $a_0$. We may suppose that
$<b,\infty>\cap W$ is contained in one component of $W\sm<a,b>$ since one side of $K$ is accessible from $\complex\sm T(X)^-_\delta$ and $a_0\in X$. If $a\in X$, let $p=a$. If not, then there exists a $\kp$ chord
$\rh$ such that $a\in\rh$. Then $\rh\subset K$ by Proposition~\ref{boundary}. Let $p$ be the endpoint of $\rh$ such that $p<a$.

Since $X\subset \ol{<x,\infty>}$ there are components of $<b,\infty>\cap W$ which are arbitrarily close to $a_0$. Choose $b<c<d$ in $K$ so that the $<c,d>$ is the closure of a component of
$W\cap<b,\infty>$ such that:
\begin{enumerate}
\item $a$ and $d$ lie in the same component of $\partial W\sm\{b,c\}$.
\item There exists $z\in<c,d>\cap X\cap W$ and an arc $I\subset \{a_o,z\}\cup [W\sm<p,d>]$ joining $a_0$ to $z$.
\item There is a $\kp$ chord $\rg\subset W$ with $z$ and $y$ as endpoints which crosses $R_t$ essentially.
Hence, $\var(f,\rg)>0$.
\item $\dm(f(\rg))<d(J^+_{a_0}\sm W,J^i_{a_0}\sm W)$.
\end{enumerate}

Conditions (1) and (2) follow because $J_{a_0}$ is a connected and closed set from $a_0$ to $\infty$ in $\{a_0\}\cup[\complex\sm T(X)_\delta^-]$ and the ray $<b,\infty>$ approaches both $a_0$ and $p$.
Conditions (3) and (4) follow from Proposition~\ref{smallarc}.
If $d\in X$, put $q=d$. Otherwise, let $q\in<d,\infty>$ such that there is a $\kp$ chord $\rh\subset K$ containing $d$ with endpoint $q$ and such that $d<q$

By a simple extension of Corollary~\ref{bumpingscc}, there exists   a bumping arc $A'$ of $T(X)$
from $p$ to $q$ such that variation is defined on each component of $A'\sm X$, $S'=A'\cup <p,q>$
is a simple closed curve with
$T(X)\subset T(S')$ and $f$ is fixed point free on $T(S')$. Since $\ol{\rg}\cap X=\{z,y\}$,
we may assume that $A'\cap \ol{\rg}=\{y\}$. Let $C$ be the arc in $\partial W$ from $a$ to $d$
disjoint from $b$. The arc $A'$ may enter $W$ and intersect $I$ several times. However, in this
case $A'$ must enter $W$ through $C$.  Since we  want to apply the Lollipop lemma,
 we will modify the arc $A'$ to a new arc $A$ which is disjoint from $I$.

Let $A$ be the set of points in $A'\cup C$ accessible from $\infty$ in $\complex\sm [S'\cup C]$.
Then $A$ is a bumping arc from $p$ to $q$, $A\cap I=\0$,
$\var(f,A)$ is defined,  $S=A\cup <p,q>$ is a simple closed curve with
$T(X)\subset T(S)$ and $f$ is fixed point free on $T(S)$. Note that $y\in A$. Then the Lollipop lemma applies
to $S$ with $R=T(<a_0,z>\cup I)$ and $R=T(I\cup <z,q>\cup A\cup <p,a_0>)$.

\smallskip
Claim: $f(z)\in R$.  Hence by Corollary~\ref{corlol}, $<a_0,z>$ contains a chord $\rg_1$ with $\var(f,\rg_1)<0$.

\smallskip
 \noindent
 Proof of Claim.
  Note that the positive direction along $\rg$ is from $z$ to $y$.  Since $z,y\in X$, $\{f(z),f(y)\}\subset X\subset T(S)=T(R)\cup T(L)$.
Choose a junction $J_y$ such that $J_{a_0}\sm W\subset J_y$ and $J_y$ runs close to $R_t$ on its way to $\rg$.  Since $\rg$ crosses $R_t$ essentially, $\var(f,\rg)>0$.
Since variation is invariant under suitable homotopies, we may assume that for each $*\in\{+,i,-\}$,
$J_y^*\sm W$ consists of exactly two components one of which is contained in $J_{a_0}$.
Let $C^*_y$ be the component of $J^*_y\sm W$ which is disjoint from $J^*_{a_0}$.
Then $C^i_y$ separates $R\cup C^+_y$ from $L\cup C^-_y$ in $\complex\sm W$
(see figure~\ref{outpic}).  Since $f(\rg)\cap J_{a_0}=\0$, if $f(z)\not\in T(R)$, $\var(f,\rg)\le 0$, a contradiction.
Hence $f(z)\in T(R)$ (and, in fact, $f(y)\in T(L)$) as desired.

Since  $f(z)\in R$,
  $<a_0,z>$ contains a chord $\rg$ with $\var(f,\rg)<0$. Repeating the same argument,
  replacing $a_0$ by $y$ and $J_0$ by $J_y$ we obtain a second chord $\rg_2$ contained in $<y,\infty>$
  such that $\var(f,\rg_2>)<0$.

  We will now show that the existence of two distinct chords $\rg_1$ and $\rg_2$ in $K$
with variation $<0$ on each leads to a contradiction.
Recall that $a_0\in\ol{<b,\infty>}$. Hence we can find $y'\in\,<b,\infty>$ with $y'\in X$ such that
$\rg_1\cup\rg_2\subset\  <a_0,y'>$ and there exists a small arc $I'\subset W'$  such that
$I'\cap <a_0,y'>=\{a_0,y'\}$.  Since $f(I')\cap J_{a_0}=\0$, $\var(f,I')=0$. We may also assume that $f$ is fixed point free
on $T(S')$, where $S'=I'\cup <a_0,y'>$.
Since $<a_0,y'>$ contains both $\rg_1$, $\rg_2$ and no chords of positive variation,
$\var(f,<a_0,y'>)\le -2$ and $\var(f,S')\le -2$. Then $\ind(f,S')=\var(f,S')+1\le -1$
a contradiction with Theorem~\ref{fpthm}. Hence $X$ has no positive geometric outchannel.

By Theorems~\ref{densechannel} and ~\ref{I=V+1}, $X$ has exactly one
negative outchannel and its variation is $-1$.

\end{proof}

Note that the following Theorem  follows from Lemma~\ref{smallarc} and Theorem~\ref{outchannel}.
\begin{thm} Suppose that $X$ is a minimal counter example to the plane  fixed point problem.
Then there exists $\delta>0$ such that the continuum $Y=T(X)_\delta^+$ is a non-separating continuum, $f$ is fixed point free on $Y$ and all accessible points of $Y$  are contained in one arc component $K$ of the boundary of $Y$.  In other words, $Y$ is homeomorphic to a disk with exactly one  channel removed which corresponds to the unique geometric outchannel of variation $-1$ of $X$. This channel compactifies on $X$.  The sides of this channel
are halflines  consisting entirely of chords of zero variation and  points in $X$. There exist arbitrarily small homeomorphisms of the tails of these halflines to the tails of $R_t$ which is the external ray corresponding to this channel.
\end{thm}

\section{Oriented maps}\label{soriented}  A {\em perfect map} \index{map!perfect}is a closed
continuous function each of whose point inverses is compact.
\emph{We will assume in the remaining sections that all maps of the plane are perfect.}
 Let $X$ and $Y$ be spaces. A
map $f:X\to Y$ is {\em monotone}  \index{map!monotone}provided for each continuum
$K\subset Y$, $f^{-1}(K)$ is connected. A map $f:X\to Y$ is {\em
confluent} \index{map!confluent} provided for each continuum $K\subset Y$ and each
component $C$ of $f^{-1}(K)$, $f(C)=K$ and   $f$ is {\em
light} \index{map!light} provided for each point $y\in Y$, $f^{-1}(y)$ is totally
disconnected.

It is well know that each homeomorphism of the plane is either orientation-preserving or
orientation-reversing. In this section we will establish an appropriate extension of
this result for confluent perfect mappings of the plane (Theorem~\ref{orient}) by
showing that such maps  either preserve or reverse local orientation. As a consequence
it follows that all perfect and confluent maps of the plane satisfy the Maximum Modulus
Theorem. We will call such maps positively- or negatively-oriented maps, respectively.
For perfect mappings of the plane, Lelek and Read have shown that confluent is
equivalent to the composition of open and monotone maps \cite{leleread74}. Holomorphic
maps are prototypes of positively-oriented maps  but positively-oriented maps, unlike
holomorphic maps, do not have to be light.  A non-separating plane continuum is said to
be {\em acyclic}. \index{acyclic}

\begin{defn}[Degree of $f_p$] Let $f:U \to \Complex$ be a map from a simply connected
domain $U\subset\complex$ into the plane.  Let $S\subset\complex$
be a positively oriented simple closed curve in $U$, and $p\in
U\setminus f^{-1}(f(S))$ a point.  Define $f_p:S\to\ucirc$ by \[
f_p(x)=\frac{f(x)-f(p)}{|f(x)-f(p)|}.\] Then $f_p$ has a well-defined {\em degree},
denoted $\degree(f_p)$. \index{degreea@$\dg(f_p)$}Note that $\degree(f_p)$ is  the  winding number
$\win(f,S,f(p))$ of $f|_S$ about $f(p)$.
\end{defn}

\begin{defn}
A map $f:U \to \Complex$ from a simply connected domain $U$ is
\emph{positively-oriented} (respectively, {\em
negatively-oriented}) provided for each simple closed curve $S$ in
$U$ and each point $p\in T(S)\setminus f^{-1}(f(S))$,
$\degree(f_p)> 0$ ($\degree(f_p)< 0$, respectively). \end{defn}
\index{map!positively oriented}\index{map!negatively oriented}
\begin{defn} A perfect surjection $f:\Complex\to\Complex$ is
\emph{oriented} provided for each simple closed curve $S$ and each
$x\in T(S)$, $f(x)\in T(f(S))$.  \index{map!oriented}
\end{defn}

Clearly every  positively oriented and  each negatively oriented map is oriented.
It will follow that all oriented maps satisfy the Maximum Modulus
Theorem~\ref{orient}.

It is well-known that both open maps and monotone maps (and hence
compositions of such maps) of continua are confluent. It will
follow (Lemma~\ref{cacyclic}) from a result of Lelek and Read
\cite{leleread74} that each perfect, oriented surjection of the plane is the
composition of a monotone map and a light open map. The following
Lemmas are in preparation for the proof of Theorem~\ref{orient}.

  \begin{lem}\label{lorient}
  Suppose $f:\Complex\to\Complex$ is a perfect surjection. Then $f$
  is confluent if and only if $f$ is oriented.
  \end{lem}

  \begin{proof}
Suppose that $f$ is oriented. Let $A$ be an arc in $\Complex$ and
let $C$ be a component of $f^{-1}(A)$. Suppose that $f(C)\not=A$.
Let $a\in A\setminus f(C)$. Since $f(C)$ does not separate $a$
from infinity, we can choose a simple closed curve $S$ with
$C\subset T(S)$, $S\cap f^{-1}(A)=\0$ and $f(S)$  so close to
$f(C)$ that $f(S)$ does not separate $a$ from $\infty$.  Then
$a\not\in T(f(S))$. Since $f$ is oriented, $f(C)\subset T(f(S))$.
Hence there exists a $y\in A\cap f(S)$.
This contradicts the fact that $A\cap f(S)=\0$. Thus $f(C)=A$.

Now suppose that $K$ is an arbitrary continuum in $\Complex$ and
let $L$ be a component of $f^{-1}(K)$. Let $x\in L$ and let $A_i$
be a sequence of arcs in $\Complex$ such that $\lim A_i=K$ and
$f(x)\in A_i$ for each $i$. Let $M_i$ be the component of
$f^{-1}(A_i)$ containing the point $x$. By the previous paragraph
$f(M_i)=A_i$. Since $f$ is perfect, $M=\limsup M_i\subset L$ is a
continuum and $f(M)=K$. Hence $f$ is confluent.

 Suppose next that $f:\Complex\to\Complex$ is not oriented. Then there exists
 a simple closed curve $S$ in $\Complex$ and
 $p\in T(S)\setminus f^{-1}(f(S))$ such that  $f(p)\not\in T(f(S))$. Let
 $L$ be a half-line with endpoint  $f(p)$ running to infinity in $\Complex\setminus
 f(S)$. Let $L^*$ be an arc in $L$ with endpoint $f(p)$ and diameter
 greater than the diameter of the continuum $f(T(S))$. Let $K$ be
 the component of $f^{-1}(L^*)$ which contains $p$. Then $K\subset
 T(S)$, since $p\in T(S)$ and $L\cap f(S)=\0$. Hence,
 $f(K)\not=L^*$, and so $f$ is not confluent.   \end{proof}

  \begin{lem}\label{finite} Let $f:\Complex\to\Complex$ be a light, open, perfect
  surjection. Then there exists an integer $k$ and a finite subset
  $B\subset\Complex$ such that $f$ is a local homeomorphism at each point of
  $\Complex\setminus B$, and for each point $y\in\Complex\setminus f(B)$,
   $|f^{-1}(y)|=k$.
   \end{lem}

\begin{proof} Let $\sphere$ be the one point compactification of
$\Complex$.  Since $f$ is perfect, we can extend $f$ to a map of $\sphere$ onto
$\sphere$  so that $f^{-1}(\infty)=\infty$. By abuse of notation we also denote the
extended map by $f$. Then $f$ is a light open mapping of the compact $2$-manifold
$\sphere$. The result now follows from a theorem of Whyburn \cite[X.6.3]{whyb42}.
\end{proof}

The following is a special case, for oriented perfect maps, of
the monotone-light factorization theorem.

\begin{lem}\label{cacyclic} Suppose that $f:\Complex\to\Complex$ is an oriented,
perfect map. It follows that $f=g\circ h$, where
$h:\Complex\to\Complex$ is a monotone perfect surjection with acyclic
fibers and $g:\Complex\to\Complex$ is a light, open perfect surjection.
\end{lem}

\begin{proof}
By the monotone-light factorization theorem \cite[Theorem
13.3]{nadl92}, $f=g\circ h$, where $h:\Complex\to X$ is monotone,
$g:X\to\Complex$ is light, and $X$ is the quotient space obtained
from $\Complex$ by identifying each component of $f^{-1}(y)$ to a
point for each $y\in\Complex$. Let $y\in \Complex$ and let $C$ be
a component of $f^{-1}(y)$. If $C$ were to separate $\Complex$,
then $f(C)=y$ would be a point while $f(T(C))$ would be a
non-degenerate continuum. Choose an arc $A\subset
\Complex\setminus \{y\}$ which meets both $f(T(C))$ and its
complement and let $x\in T(C)\setminus C$ such that $f(x)\in A$.
If $K$ is the component of $f^{-1}(A)$ which contains $x$, then
$K\subset T(C)$. Hence $f(K)$ cannot map onto $A$ contradicting
the fact that $f$ is confluent. Thus, for each $y\in\Complex$, each
component of $f^{-1}(y)$ is acyclic.

By Moore's Plane Decomposition Theorem \cite{dave86}, $X$ is
homeomorphic to $\Complex$. Since $f$ is confluent, it is easy to
see that $g$ is confluent. By a theorem of Lelek and
Read~\cite{leleread74} $g$ is open since it is confluent and light
(also see \cite[Theorem 13.26]{nadl92}).  Since $h$ and $g$ factor
the perfect map $f$ through a Hausdorff space $\Complex$, both $h$
and $g$ are perfect \cite[3.7.5]{enge89}.
\end{proof}

\begin{thm}[Maximum Modulus Theorem] \label{orient}\index{Maximum Modulus Theorem}
Suppose that $f:\Complex\to\Complex$ is a perfect surjection. Then the following are
equivalent:
\begin{enumerate}
\item\label{pnorient} $f$ is either  positively or
 negatively oriented. \item \label{iorient}$f$ is
oriented. \item\label{conf}  $f$ is confluent.
\end{enumerate} Moreover, if $f$ is oriented, then for any non-separating continuum $X$,
$\bd(f(X))\subset f(\bd(X))$.
\end{thm}

\begin{proof}
It is clear that (\ref{pnorient}) implies (\ref{iorient}). By
Lemma~\ref{lorient} every oriented map is confluent. Hence suppose
that $f:\Complex\to\Complex$ is a perfect confluent map. By
Lemma~\ref{cacyclic}, $f=g\circ h$, where $h:\Complex\to\Complex$
is a monotone perfect  surjection with acyclic fibers and
$g:\Complex\to\Complex$ is a light, open perfect surjection. By Stoilow's
Theorem \cite{whyb64} there exists a homeomorphism
$j:\Complex\to\Complex$ such that $g\circ j$ is an analytic surjection.
Then $f=g\circ h=(g\circ j)\circ (j^{-1}\circ h)$. Since
$k=j^{-1}\circ h$ is a monotone surjection of $\Complex$ with
acyclic fibers, it is a near homeomorphism \cite[Theorem 25.1]{dave86}. That is, there exists
a sequence $k_i$ of homeomorphisms of $\Complex$ such that $\lim
k_i=k$. We may assume that all of the
$k_i$ have the same orientation.

Let $f_i=(g\circ j)\circ k_i$,
$S$ a simple closed curve in the domain of $f$ and $p\in T(S)\setminus f^{-1}(f(S))$.
Note that $\lim f^{-1}_i(f_i(S))\subset f^{-1}(f(S))$. Hence $p\in
T(S)\setminus f^{-1}_i(f_i(S))$ for $i$ sufficiently large.
Moreover, since $f_i$ converges to $f$, $f_i|_S$ is homotopic to
$f|_S$ in the complement of $f(p)$ for $i$ large. Thus for large
$i$, $\degree((f_i)_p)=\degree(f_p)$, where
\[(f_i)_p(x)=\frac{f_i(x)-f_i(p)}{|f_i(x)-f_i(p)|}
\text{ and } f_p(x)=\frac{f(x)-f(p)}{|f(x)-f(p)|}.\] Since $g\circ
j$ is an analytic map, it is positively oriented and
$\degree((f_i)_p)=\degree(f_p)>0$ if $k_i$ is orientation
preserving and $\degree((f_i)_p)=\degree(f_p)<0$ if $k_i$ is
orientation reversing. Thus, $f$ is positively-oriented if each
$k_i$ is orientation-preserving and $f$ is negatively-oriented if
each $k_i$ is orientation-reversing.

Suppose that $X$ is a non-separating continuum and $f$ is oriented. Let
$y\in\bd(f(X))$. Choose $y_i\in\bd(f(X)$ and rays $R_i$, joining $y_i$ to $\infty$
such that $R_i\cap f(X)=\{y_i\}$ and $\lim y_i=y$. Choose $x_i\in X$ such that
$f(x_i)=y_i$. Since $f$ is confluent, there exists
closed and connected sets $C_i$, joining $x_i$ to $\infty$ such that $C_i\cap X\subset f^{-1}(y_i)$.
Hence there exist  $x'_i\in f^{-1}(y_i)\cap \bd(X)$. We may assume that $\lim x'_i=x_\infty\in\bd(X)$
and $f(x_\infty)=y$ as desired.
\end{proof}

 We shall need the following three results in the
next section.

\begin{lem}\label{acyclic} Let $X$ be a plane continuum and
$f:\Complex\to\Complex$ a perfect, surjective map such that $f^{-1}(T(X))= T(X)$ and
$f|_{\Complex\setminus  f^{-1}(T(X))}$ is confluent.  Then for each
$y\in\Complex\setminus T(X)$,  each component of $f^{-1}(y)$ is
acyclic.
\end{lem}
\begin{proof}
Suppose there exists $y\in\Complex\setminus T(X)$ such that some
component $C$ of $f^{-1}(y)$ is not acyclic. Then there exists
$z\in T(C)\setminus f^{-1}(y)\cup T(X)$. By unicoherence of
$\Complex$, $T(X)\cup \{y\}$ does not separate $f(z)$ from infinity
in $\Complex$. Let $L$ be a ray in $\Complex\setminus[T(X)\cup
\{y\}]$ from $f(z)$ to infinity. Then $L=\cup L_i$, where each
$L_i\subset L$ is an arc with endpoint $f(z)$. For each $i$ the
component $M_i$ of $f^{-1}(L_i)$ containing $z$ maps onto $L_i$.
Then $M=\cup M_i$ is a connected closed subset in
$\Complex\setminus f^{-1}(y)$ from $z$ to infinity. This is a
contradiction since $z$ is contained in a bounded complementary component of $f^{-1}(y)$.
\end{proof}

\begin{thm}\label{confeq}
Let $X$ a plane continuum and
$f:\Complex\to\Complex$ a perfect, surjective map such that $f^{-1}(T(X))=T(X)$ and
$f|_{\Complex\setminus f^{-1}(T(X))}$ is confluent.
 If $A$ and $B$ are crosscuts of $T(X)$ such that
 $B\cup X$ separates $A$ from $\infty$ in $\Complex$, then
  $f(B)\cup T(X)$ separates $f(A)\setminus f(B)$ from $\infty$.
\end{thm}
\begin{proof} Suppose not. Then there exists a half-line  $L$
joining $f(A)$ to infinity in $\Complex\setminus (f(B)\cup T(X))$. As
in the proof of Lemma~\ref{acyclic}, there exists a closed and
connected set $M\subset \Complex\setminus (B\cup X)$ joining $A$
to infinity, a contradiction.
\end{proof}

\begin{prop}\label{ray} Under the conditions of
Theorem~\ref{confeq}, if $L$ is a ray irreducible from $T(X)$ to
infinity, then each component of $f^{-1}(L)$  is closed in $\complex\sm X$
and is a  connected set from $X$ to
infinity.
\end{prop}

\section{Induced maps of prime ends} Suppose that
$f:\Complex\to\Complex$ is an oriented perfect surjection and
$f^{-1}(Y)=X$, where $X$ and $Y$ are acyclic continua.    We will
show that in this case the map $f$ induces a confluent map $F$ of
the circle of prime ends of $X$ to the circle of prime ends of
$Y$. This result was announced by Mayer in the early 1980's but
never appeared in print. It was also used (for homeomorphisms) by
Cartwright and Littlewood in \cite{cartlitt51}.  There are easy
counterexamples that show if $f$ is not confluent then it may not
induce a continuous function between the circles of prime ends.

\begin{thm}\label{tinduced}\index{map!on circle of prime ends}
Let $X$ and $Y$ be non-degenerate acyclic plane continua and
$f:\Complex\to\Complex$ a perfect map such that:
\begin{enumerate}
\item \label{cut} $Y$ has no cut point, \item \label{complement}
$f^{-1}(Y)=X$ and \item \label{sep}  $f|_{\Complex\setminus X}$ is
confluent.
\end{enumerate}
Let $\varphi: \disk^\infty\to \sphere\setminus X$ and
$\psi:\disk^\infty\to \sphere\setminus Y $ be conformal mappings. Define
$\hat{f}:\disk^\infty\to\disk^\infty$ by $\hat{f}=\psi^{-1} \circ f \circ
\varphi$.

Then $\hat{f}$ extends to a map $\bar{f}:\ol{\disk^\infty}\to\ol{\disk^\infty}$.
Moreover, $\bar{f}^{-1}(S^1)=S^1$ and $F=\bar{f}|_{S^{1}}$ is a
confluent map.
\end{thm}

\begin{proof} Note that $f$ takes accessible points of $X$ to
accessible points of $Y$. For if $P$ is a path in
$[\Complex\setminus X]\cup \{p\}$ with endpoint $p\in X$, then by
(\ref{complement}), $f(P)$ is a path in $[\Complex\setminus Y]\cup
\{f(p)\}$ with endpoint $f(p)\in Y$.

Let $A$ be a crosscut of $X$ such that the diameter of $f(A)$ is less than half of the
diameter of $Y$ and let $U$ be the bounded component of $\Complex\setminus (X\cup A)$.
Let the  endpoints of $A$ be $x,y\in X$ and suppose that $f(x)=f(y)$. If  $x$ and $y$
lie in the same component of $f^{-1}(f(x))$   then each crosscut $B\subset U$ of $X$ is
mapped to a  generalized return cut of $Y$ based at $f(x)$ (i.e., the endpoints of $B$
map to $f(x)$) by (1). Note that in this case by (1), $\bd{f(U)}\subset
f(A)\cup \set{f(x)}$.

 Now suppose that $f(x)=f(y)$ and $x$
and $y$ lie in distinct components of $f^{-1}(f(x))$. Then by unicoherence of
$\Complex$, $\bd{U}\subset A\cup X$ is a connected set and $\bd{U}\not\subset
\bar{A}\cup f^{-1}(f(x))$. Now $\bd{U}\setminus (\bar{A}\cup
f^{-1}(f(x)))=\bd{U}\setminus f^{-1}(f(\bar{A}))$ is an open non-empty set in $\bd{U}$ by(2).  Thus
there is a crosscut $B\subset U\setminus f^{-1}(f(\bar{A}))$ of $X$ with
$\bar{B}\setminus B\subset\bd{U}\setminus f^{-1}(f(\bar{A}))$. Now $f(B)$ is contained
in a bounded component of $\Complex\setminus (Y\cup f(A))=\Complex\setminus (Y\cup
f(\bar{A}))$ by Theorem~\ref{confeq}. Since $Y\cap f(\bar{A})=\{f(x)\}$ is connected and
$Y$ does not separate $\Complex$, it follows by unicoherence that $f(B)$ lies in a
bounded component of $\Complex\setminus f(\bar{A})$. Since $Y\setminus \{f(x)\}$ meets
$f(\bar{B})$ and misses $f(\bar{A})$ and $Y\setminus \{f(x)\}$ is connected,
$Y\setminus\{f(x)\}$ lies in a bounded complementary component of $f(\bar{A})$. This is
impossible as the diameter of $f(A)$ is smaller than the diameter of $Y$.  It follows
that there exists a $\delta>0$ such that if the diameter of $A$ is less than $\delta$
and $f(x)=f(y)$, then $x$ and $y$ must lie in the same component of $f^{-1}(f(x))$.

  In order to define the extension $\bar{f}$ of $f$ over the boundary $S^1$ of
  $\disk^\infty$, let $C_i$ be a chain  of crosscuts of $\disk^\infty$ which
  converge to a point $p\in S^1$ such that $A_i=\varphi(C_i)$
  is a null chain of crosscuts or return cuts of $X$ with endpoints $a_i$ and
  $b_i$ which converge to a point $x\in X$. There are three cases
  to consider:

  Case 1. $f$ identifies the endpoints of $A_i$ for some $A_i$ with
  diameter less than $\delta$. In this case the chain of crosscuts
  is mapped by $f$ to a chain of generalized return cuts based at
  $f(a_i)=f(b_i)$. Hence $f(a_i)$ is an accessible point of $Y$
  which corresponds (under $\psi^{-1}$) to a unique point $q\in S^1$ (since $Y$ has no cutpoints).
  Define $\bar{f}(p)=q$.

Case 2. Case 1 does not apply and there exists an infinite subsequence $A_{i_{j}}$ of
crosscuts such that
  $f(\bar{A}_{i_{j}})\cap f(\bar{A}_{i_{k}})=\0$ for $j\not=k$. In this case
 $f(A_{i{_j}})$ is  a chain of generalized crosscuts which
converges to
  a point $f(x)\in Y$. This chain corresponds to a unique point
  $q\in S^1$ since $Y$ has no cut points. Define $\bar{f}(p)=q$.

  Case 3. Cases 1 and 2 do not apply. Without loss of generality suppose
   there exists an $i$ such that
  for $j>i$ $f(\bar{A}_i)\cap
  f(\bar{A}_j)$ contains  $f(a_i)$.
  In this case $f(A_j)$ is a chain of generalized
  crosscuts based at the accessible point $f(a_i)$ which
  corresponds to a unique point $q$ on $S^1$ as above. Define
  $\bar{f}(p)=q$.

  It remains to be shown that $\bar{f}$ is a continuous extension of
  $\hat{f}$ and  $F$ is confluent. For continuity it suffices to show
  continuity at $S^1$. Let $p\in S^1$ and let $C$ be a small
  crosscut of $\disk^{\infty}$ whose endpoints are on opposite sides of $p$ in $S^1$ such that
  $A=\varphi(C)$ has diameter less than $\delta$ \cite{miln00} and such
  that the endpoints of $A$ are two accessible points of $X$.
  Since $f$ is uniformly continuous near $X$, the diameter of
  $f(A)$ is small and since $\psi^{-1}$ is uniformly continuous with
  respect to connected sets in the complement of $Y$
  (\cite{urseyoun51}), the diameter of
  $B=\psi^{-1} \circ f\circ \varphi$ is small. Also $B$ is
  either a generalized crosscut or generalized return cut.
  Since $\hat{f}$ preserves separation of crosscuts,
  it follows that the image of
  the domain $U$ bounded by $C$ which does not contain $\infty$ is
  small. This implies continuity of $\bar{f}$ at $p$.

  To see that $F$ is confluent let $K\subset S^1$ be a
  subcontinuum and let $H$ be a component of $\bar{f}^{-1}(K)$.
  Choose a chain  of crosscuts $C_i$ such that $\varphi(C_i)=A_i$ is a crosscut of $X$
  meeting $X$ in two accessible points $a_i$ and $b_i$, $C_i\cap
  \bar{f}^{-1}(K)=\0$ and $\lim C_i=H$. It follows from the preservation
  of crosscuts (see Theorem~\ref{confeq}) that $\hat{f}(C_i)$ separates $K$ from $\infty$.
  Hence $\hat{f}(C_i)$ must meet $S^1$ on both sides of $K$ and
  $\lim \bar{f}(C_i)=K$. Hence $F(H)=\lim \bar{f}(C_i)=K$ as required.
\end{proof}

\begin{cor} Suppose that $f:\Complex\to\Complex$ is a perfect, oriented  map
 of the plane,   $X\subset\Complex$ is  a subcontinuum
without cut points and $f(X)=X$. Let  $\hX$ be the component of
$f^{-1}(f(X))$ containing $X$. Let $\varphi: \disk^\infty \to\sphere\setminus
T(\hX) $ and $\psi:\disk^\infty\to\sphere\setminus T(X)$ be
conformal mappings. Define $\hat{f}:\disk^\infty\sm
\varphi^{-1}(f^{-1}(X))\to\disk^\infty$ by $\hat{f}=\psi^{-1} \circ f \circ
\varphi$. Put $S^1=\partial\disk^\infty$.

Then $\hat{f}$ extends over $S^1$ to a map
$\bar{f}:\ol{\disk^\infty}\to\ol{\disk^\infty}$. Moreover $\bar{f}^{-1}(S^1)=S^1$ and
$F=\bar{f}|_{S^{1}}$ is a confluent map.
\end{cor}

\begin{proof}
By Lemma~\ref{cacyclic} $f=g\circ m$ where $m$ is a monotone
perfect and onto mapping of the plane with acyclic point inverses,
and $g$ is an open and perfect surjection of the plane to itself.
By Lemma~\ref{finite}, $f^{-1}(X)$ has finitely many components.
It follows that
there exist a simply connected open set $V$, containing $T(X)$, such that if $U$ is the component
of $f^{-1}(V)$ containing $\hX$, then $U$ contains no other components of $f^{-1}(X)$.
It is easy to see that $f(U)=V$ and that $U$ is simply connected.
Hence $U$ and $V$ are homeomorphic
 to $\Complex$. Then  $f|_U:U\to V$ is a  confluent map.
 The result now follows from Theorem~\ref{tinduced}
 applied to $f$ restricted to $U$.
\end{proof}

\section{Fixed points for positively oriented maps}
In this section we will consider a positively oriented map of the
plane. As we shall see below, a  straight forward application of
the tools developed above will give us the desired fixed point
result.
 We will assume, by way of contradiction, that $f:\Complex \to \Complex$ is
 a positively oriented map, $X$ is a  plane
 continuum such that $f(X)\subset T(X)$ and $T(X)$ contains no fixed points of $f$.

\begin{lem}\label{closed}  Let $f:\Complex\to\Complex$ be a map and $X\subset\complex$ a
 continuum such that $f(X)\subset T(X)$. Suppose $C=(a,b)$ is a crosscut of
the continuum $T(X)$. Let  $v\in (a,b)$ be a point and $J_v$ be a junction such that
$J_v\cap (X\cup C)=\{v\}$. Then there exists an arc $I$ such that $S=I\cup C$ is a
simple closed curve, $T(X)\subset T(S)$ and $f(I)\cap J_v= \0$.
\end{lem}

\begin{proof}
Since $f(X)\subset T(X)$ and $J_v\cap X=0$, it is clear that there exists an arc $I$ with
endpoints $a$ and $b$ sufficiently close to  $T(X)$ such that $I\cup C$ is a simple closed
curve, $T(X)\subset T(I\cup C)$ and $f(I)\cap J_v=\0$. This completes the proof.
\end{proof}

\begin{cor}\label{posvar}
Suppose $X\subset\complex$ is a continuum, $f:\Complex\to \Complex$
 a positively oriented
map such that $f(X)\subset T(X)$. Then for each crosscut
  $C$  of $T(X)$ such that $f(\cl C)\cap \cl C=\0$, $\var(f,C)\geq 0$
\end{cor}

\begin{proof}
Suppose that $C=(a,b)$ is a crosscut of $T(X)$ such that $f(\cl C)\cap \cl C=\0$ and
$\var(f,C)\not=0$. Choose a junction $J_v$ such that $J_v\cap (X\cup C)=\set{v}$ and
$v\in C\setminus X$. By Lemma~\ref{closed}, there exists an arc $I$ such that $S=I\cup
C$ is a simple closed curve and $f(I)\cap J_v=\0$. Moreover, by choosing $I$
sufficiently close to $X$, we may assume that $v\in\Complex\setminus f(S)$. Hence
$\var(f,C)=\text{Win}(f,S,v)\geq 0$ by the remark following Definition~\ref{vararc}.
\end{proof}

\begin{thm}\label{fixpoint}\index{fixed point!for positively oriented maps}
Suppose $f:\Complex\to\Complex$ is a positively oriented map and
$X$ is a  continuum such that $f(X)\subset T(X)$. Then
there exists a point $x_0\in T(X)$ such that $f(x_0)=x_0$.
\end{thm}

\begin{proof}
Suppose  we are given a  continuum
 $X$ and $f:\Complex\to \Complex$  a positively oriented map such that
 $f(X)\subset T(X)$.
  Assume that $f|_{T(X)}$ is fixed point free.
 Choose a simple closed curve $S$ such that $X\subset T(S)$
 and points $a_0<a_1<\ldots<a_n$
 in $S\cap X$
  such that for each $i$ $C_i=(a_i,a_{i+1})$ is a sufficiently small crosscut of
 $X$,
 $f(\ol{C_i})\cap \ol{C_i}=\0$ and  $f|_{{T(S)}}$ is fixed point free. By Corollary~\ref{posvar},
 $\var(f,C_i)\geq 0$ for each $i$. Hence by Theorem~\ref{I=V+1},
 $\ind(f,S)=\sum \var(f,C_i) +1 \geq 1$. This contradiction with Theorem~\ref{fpthm}
 completes the proof.
 \end{proof}

 \begin{cor} \index{point of period two!for oriented maps}
 Suppose $f:\Complex\to\Complex$  is a perfect, oriented
 map
 and $X$ is a continuum such that
 $f(X)\subset T(X)$. Then
there exists a point $x_0\in  T(X)$ of period at most  2.
\end{cor}
 \begin{proof}
By Theorem~\ref{orient}, $f$ is either positively or negatively
oriented. In either case, the second iterate $f^2$ is positively
oriented and must have a fixed point in $T(X)$ by
Theorem~\ref{fixpoint}.
\end{proof}

\nocite{bell76} \nocite{bell78}
\bibliographystyle{plain}
\bibliography{c:/lex/references/refshort}

\printindex
\begin{theindex}

  \item accessible point, 14
  \item acyclic, 32
  \item allowable partition, 13

  \indexspace

  \item $\B$, 15
  \item $\B^\infty$, 15
  \item bumping
    \subitem arc, 11
    \subitem simple closed curve, 11

  \indexspace

  \item $C(a,b)$, 18
  \item Carath\'eodory Loop, 12
  \item chain of crosscuts, 13
    \subitem equivalent, 13
  \item channel, 14
    \subitem dense, 14, 26
  \item completing a bumping arc, 11
  \item $\complex$, 2
  \item $\rsphere$, 2
  \item $\ECH(K)$, 15
  \item $\HCH(B\cap K)$, 15
  \item convex hull
    \subitem Euclidean, 15
    \subitem hyperbolic, 15
  \item counterclockwise order
    \subitem on an arc in a simple closed curve, 3
  \item crosscut, 11
    \subitem shadow, 12

  \indexspace

  \item $\Disk$, 18
  \item defines variation near $X$, 24
  \item degree, 3
  \item $\dg(g)$, 3
  \item $\dg(f_p)$, 32
  \item $\partial$ boundary operator, 2
  \item $\disk^\infty$, 13

  \indexspace

  \item embedding
    \subitem orientation preserving, 4
  \item $\mc{E}_t$, 13
  \item external ray, 14
    \subitem end of, 14
    \subitem essential crossing, 14
    \subitem landing point, 14

  \indexspace

  \item fixed point
    \subitem for positively oriented maps, 37
  \item $(f,X,\eta)$, 24

  \indexspace

  \item $\rG$, 18
  \item $\fg$, 18
  \item $\rg$, 15
  \item $\Gamma$, 18
  \item gap, 17
  \item geodesic
    \subitem hyperbolic, 18

  \indexspace

  \item hull
    \subitem hyperbolic, 17
  \item hyperbolic
    \subitem geodesic, 15
    \subitem halfplane, 15
  \item hyperbolic chord, 18

  \indexspace

  \item $id$ identity map, 3
  \item $\im(\mc{E}_t)$, 14
  \item impression, 14
  \item index, 3
    \subitem fractional, 4
    \subitem I=V+1 for Carath\'eodory Loops, 13
    \subitem Index=Variation+1 Theorem, 6
  \item $\ind(f,A)$, 4
  \item $\ind(f,g)$, 3
  \item $\ind(f,g_{[a,b]})$, 4
  \item $\ind(f,S)$, 4

  \indexspace

  \item junction, 5
  \item J{\o}rgensen Lemma, 18

  \indexspace

  \item $\KP$, 17
  \item $\KP_\delta$, 20
  \item $\kp$ chord, 17
  \item $\KP^{\pm}_\delta$, 22
  \item $\KPP$, 17
  \item $\KPP_\delta$, 20
  \item Kulkarni-Pinkall
    \subitem Lemma, 16
    \subitem Partition, 17

  \indexspace

  \item landing point, 14
  \item Lollipop Lemma, 9

  \indexspace

  \item map
    \subitem confluent, 31
    \subitem light, 31
    \subitem monotone, 31
    \subitem negatively oriented, 32
    \subitem on circle of prime ends, 34
    \subitem oriented, 32
    \subitem perfect, 31
    \subitem positively oriented, 32
  \item maximal ball, 15
  \item Maximum Modulus Theorem, 33

  \indexspace

  \item narrow strip, 27

  \indexspace

  \item order
    \subitem on subarc of simple closed curve, 3
  \item orientation preserving
    \subitem embedding, 4
  \item outchannel, 24
    \subitem geometric, 24
    \subitem uniqueness, 28

  \indexspace

  \item point of period two
    \subitem for oriented maps, 37
  \item $\pr(\mc{E}_t)$, 14
  \item prime end, 13
    \subitem channel, 14
    \subitem impression, 14
    \subitem principal continuum, 14
  \item principal continuum, 14

  \indexspace

  \item $\real$, 2
  \item $R_t$, 14

  \indexspace

  \item shadow, 12
  \item $\Sh(A)$, 12
  \item smallest ball, 15
  \item standing hypothesis, 24

  \indexspace

  \item topological hull, 2
  \item $T(X)$, 2
  \item $T(X)_\delta$, 22

  \indexspace

  \item $U^\infty$, 2

  \indexspace

  \item $\var(f,A)$, 12
  \item variation
    \subitem for crosscuts, 11
    \subitem of a simple closed curve, 6
    \subitem on an arc, 5
    \subitem on finite union of arcs, 6
  \item $\var(f,A,S)$, 5
  \item $\var(f,S)$, 6

  \indexspace

  \item $\win(g,S^1,w)$, 3

\end{theindex}

\end{document}